\documentclass{amsart}
\usepackage{geometry}
\usepackage[utf8]{inputenc}
\usepackage{amsmath,amsfonts,amssymb,amsthm,tabto,commath,multicol,tikz-cd}
\usepackage{diagrams}
\usetikzlibrary{arrows,shapes,automata,backgrounds,petri,fit}
\usepackage{hyperref}

\newtheorem{thm}{Theorem}[section]
\newtheorem{lem}[thm]{Lemma}
\newtheorem{prop}[thm]{Proposition}
\newtheorem{coro}[thm]{Corollary}
\newtheorem{conj}[thm]{Conjecture}

\theoremstyle{definition}

\newtheorem{ex}[thm]{Example}

\newtheorem{rem}[thm]{Remark}
\newtheorem{openq}[thm]{Question}
\newtheorem{defn}[thm]{Definition}

\newcommand*{\onto}{\ensuremath{\twoheadrightarrow}}
\newcommand*{\into}{\ensuremath{\hookrightarrow}}

\newcommand{\xto}[1]{\xrightarrow{#1}}

\DeclareMathOperator{\Hom}{\mathrm Hom}

\DeclareMathOperator{\End}{\mathrm End}
\DeclareMathOperator{\rk}{\mathrm{rk}}

\newcommand{\poset}{\mathcal{P}}
\newcommand{\incidence}{I(\poset)}
\newcommand{\X}{\mathcal{X}}
\newcommand{\Y}{\mathcal{Y}}
\newcommand{\E}{\mathcal{E}}
\newcommand{\mods}{\mathsf{mod}}

\newcommand{\add}{\mathrm{add}}
\newcommand{\undim}{\underline{\mathrm{dim}}}
\newcommand{\unrk}{\underline{\mathrm{rk}}}

\title{Homological Approximations in Persistence Theory}

\author{Benjamin Blanchette}
\address{Benjamin Blanchette, D{\'e}partment de Math{\'e}matiques,
Universit{\'e} de Sherbrooke,
Sherbrooke, QC, J1K 2R1, Canada}
\email{Benjamin.Blanchette@USherbrooke.ca}

\author{Thomas Br{\"u}stle}
\address{Thomas Br{\"u}stle, D{\'e}partment de Math{\'e}matiques,
Universit{\'e} de Sherbrooke,
Sherbrooke, QC, J1K 2R1, Canada and
 Bishop's University, Sherbrooke, QC, J1M 1Z7, Canada}
\email{Thomas.Brustle@USherbrooke.ca, tbruestl@bishops.ca}

\author{Eric J. Hanson}
\address{Eric J. Hanson, LACIM, Universit{\'e} du Qu{\'e}bec {\`a} Montr{\'e}al, Montr{\'e}al, QC, H2L 2C4, Canada}
\email{ejhanso3@ncsu.edu}

\date{12 December, 2022}


\keywords{persistence modules, invariants, Grothendieck groups, relative homological algebra, exact structures.\\
\indent This is the accepted manuscript of a paper published in {\it Canadian Journal of Mathematics}. This version is available under a Creative Commons CC BY-NC-ND license. The official version of record is available at \url{https://doi.org/10.4153/S0008414X22000657}.  \copyright \ The Author(s)}

\begin{document}

\maketitle

\begin{abstract}
    We define a class of invariants, which we call homological invariants, for persistence modules over a finite poset. Informally, a homological invariant is one that respects some homological data and takes values in the free abelian group generated by a finite set of indecomposable modules. We focus in particular on groups generated by ``spread modules'', which are sometimes called ``interval modules'' in the persistence theory literature. We show that both the dimension vector and rank invariant are equivalent to homological invariants taking values in groups generated by spread modules. We also show that the free abelian group generated by the ``single-source'' spread modules gives rise to a new invariant which is finer than the rank invariant. They are also thankful to an anonymous referee for their thorough reading of this paper and suggestions for improvement.
\end{abstract}

\tableofcontents

\section{Introduction}
When taken at face value, a finite set of points in Euclidean space has no interesting topological features. There are many cases, however, where such a set of points can be seen as approximating something more interesting. For a simple example, we can imagine a set of points in $\mathbb{R}^2$ which are randomly sampled from the unit circle. In such a context, one may wish to speak of the ``topology of the points" as the topology of this approximated space. The theory of \emph{persistent homology}, which originated in the 1990s (see \cite{frosini,robins,DE,barannikov}, and 10 years later \cite{ELZ}), gives a formalization of this concept. The core idea is to construct a series of topological spaces (often simplicial complexes) which begins as a discrete set of points and ends as a contractible space. This series of topological spaces, often called a filtration, comes equipped with inclusion maps. The homological features which ``persist'' under several subsequent inclusions are considered to be those ``of interest''.

    The homology, often with field coefficients, of these topological spaces and the induced maps between them form what is called a \emph{persistence module}. These objects can also be formulated abstractly. Motivated by works such as \cite{BS,CL,BC-B}, let $\mathcal{P}$ be a poset and $K$ a field. We consider $\mathcal{P}$ as a category, with a unique morphism $a \rightarrow b$ whenever $a \leq b \in \mathcal{P}$. A \emph{persistence module over $\mathcal{P}$} is then defined to be a functor $M$ from $\mathcal{P}$ to the category of $K$-vector spaces. We adopt the standard convention of denoting by $M(a)$ the result of applying $M$ to the object $a \in \poset$ and by $M(a,b)$ the result of applying $M$ to the unique morphism $a \rightarrow b$ in $\poset$. If $\dim M(a) < \infty$ for all $a \in \mathcal{P}$, the persistence module $M$ is said to be \emph{pointwise finite-dimensional} (pfd). All of the persistence modules considered in this paper will be assumed to have this property unless otherwise stated.
  
Let us now restrict to the case where $|\poset| < \infty$. Then the \emph{incidence algebra} of $\poset$ is the finite-dimensional $K$-algebra $\incidence:= KQ/J$, where $Q$ is the Hasse quiver of $\poset$ and $J$ is the two-sided ideal generated by all commutativity relations. Moreover, it is well-known that the category $\mods\incidence$ of finite-dimensional (right) $\incidence$-modules is isomorphic to that of pfd persistence modules over $\poset$. As such, we will freely move between these categories throughout this paper. We refer to \cite{ARS,ASS} for more details about finite-dimensional algebras and their module categories, to \cite{simson} for more details about incidence algebras, and to \cite[Chapter~4]{ec1} for more details about general poset theory.

When $\poset$ is totally ordered, the category $\mods\incidence$ is very well behaved. Indeed, in this case, the category $\mods\incidence$ contains finitely many indecomposable objects, indexed by pairs $a \leq b \in \mathcal{P}$. These are the so-called ``interval modules", see Definition~\ref{def:spreadRep} below. Every persistence module can then be uniquely written as a direct sum of these indecomposables, yielding what is often referred to as the ``barcode'' of the persistence module.  When viewed in this way, the parameters $a$ and $b$ describing the ends of a ``bar'' can be seen as the ``birth'' and ``death'' of some topological feature. For additional details on this topic, readers are referred to \cite[Chapter~1]{oudot},  \cite{C-B}, and \cite{IRT}. In Section~\ref{sec:barcode}, we will also explicitly describe how barcodes are related to the invariants introduced in this paper.

In both applications and theory, the natural next step is to ask what happens when the poset $\poset$ is no longer totally ordered. This gives rise to the theory of ``multiparameter persistence'' \cite{CZ}, in which the poset $\mathcal{P}$ is often taken to be a closed interval in the lattice $\mathbb{Z}^n$. Except in small cases, this amounts to studying algebras of ``wild representation type''; i.e., algebras for which a complete description of all indecomposable modules is impossible. Since it is no longer  computationally feasible to describe arbitrary persistence modules in terms of direct sums of indecomposables in this setting, one instead can turn to the study of \emph{invariants}. The invariants we consider are maps $p:\mods \incidence \rightarrow \mathbb{Z}^n$ which are constant on isomorphism classes. We further assume that our invariants are \emph{additive}, meaning that $p(M\oplus N) = p(M) + p(N)$. See Definition~\ref{def:invariant} for precise definitions.

Perhaps the most natural invariant to consider is the \emph{dimension vector}, also known as the \emph{Hilbert function}. Given a persistence module $M$, the \emph{dimension vector} is defined by
    $$\undim M = (\dim M(a))_{a \in \mathcal{P}}.$$
As discussed in Section~\ref{sec:motivation}, for each $a \in \poset$, there is an indecomposable projective $\incidence$-module $P_a$ which satisfies
$$\dim \Hom(P_a,M) = \dim M(a).$$
(One may also refer to $P_a$ as a ``free persistence module". See Section~\ref{sec:terminology} for notes on our choice of terminology.) More generally, given any finite set $\mathcal{X}$ of indecomposable $\incidence$-modules, one obtains an invariant $(\dim\Hom(R,M))_{R \in \mathcal{X}}$. We will call invariants of this form \emph{dim-hom invariants}. These are the types of invariants we consider in this paper.

In order to motivate the study of dim-hom invariants, we return to our discussion of classical invariants. While the dimension vector is easy to compute, it is agnostic to the maps comprising the persistence module. In persistence theory, however, we are often interested in how long topological features \emph{persist} through a filtration. By construction, this is precisely the information stored in the maps comprising the module. One way to incorporate this information is to use the \emph{rank invariant}, defined by
 \begin{equation}\label{eqn:rankInv}\unrk M = (\rk M(a,b))_{a\leq b \in \poset},\end{equation}
where $\rk M(a,b)$ denotes the rank of the linear map $M(a,b)$ in the traditional sense.
 This invariant was first introduced in \cite{CZ}, where it was the first invariant tailored for multiparameter persistence homology. We note the rank invariant is strictly finer than the dimension vector, since for all $a \in \poset$ one has $\dim M(a) = \rk M(a,a)$. (It is straightforward to construct modules with the same dimension vector but different rank invariants.)
 
 When $\poset$ is totally ordered, the rank invariant is equivalent to the barcode, see \cite[Theorem~12]{CZ}. This in particular means that the rank invariant is \emph{complete} in this case; i.e., that the isomorphism class of a persistence module can be recovered from its barcode. For more general posets, however, there exist nonisomorphic persistence modules with the same rank invariant. See e.g. Example~\ref{ex:recVrk}.

In the recent preprint \cite{BOO}, the rank invariant is shown to be a dim-hom invariant. While an interesting result on its own, the authors further give a new interpretation of the rank invariant using relative homological algebra. More precisely, let $M$ be a module. The authors show that $\unrk M$ is equivalent to the representative $[M]_\X$ of $M$ in the \emph{relative Grothendieck group} of $\mods \incidence$ with respect to some exact structure. See Section~\ref{sec:relHom} for background and definitions pertaining to relative Grothendieck groups and Theorem~\ref{thm:boo} for the precise statement proven in \cite{BOO}. In particular, the representative $[M]_\X$ can be interpreted as a ``signed approximation" of $M$ by the ``relative projective'' objects in the exact structure. In this paper, we take the reverse approach: we first develop a theoretical framework using homological algebra, and then we build new invariants using it. We refer to those invariants which fit into this framework as \emph{homological invariants}.

We show that the dimension vector, the rank invariant, and, when the Hasse quiver of $\poset$ is Dynkin type $A$, the barcode are all examples of homological invariants. We also concretely describe a new homological invariant for persistence modules over an arbitrary finite poset. Unless the poset satisfies a very restrictive condition, slightly more general than being totally ordered, this new invariant is strictly finer than the rank invariant. Finally, we compare the approximations coming from our homological invariants to other notions of approximation from the persistence theory literature. See Section~\ref{sec:motivation} for additional details and motivation.

%%%%%%%%%%%%%%%%%%%%%%%%%%%%%%%%%%%%%%%%%%%%%%%%%%%%%%%%

\subsection{Organization and main results}
The contents of this paper are as follows.

In Section~\ref{sec:terminology}, we fix notation and terminology for use in the later sections. As we are mainly drawing intuition from the representation theory of finite-dimensional algebras, this terminology sometimes diverges from that in the persistence theory literature. In particular, we define \emph{spread modules} (Definition~\ref{def:spreadRep}), which are typically referred to as ``interval modules'' in the persistence theory literature.

In Section~\ref{sec:motivation}, we provide motivation for our study of homological invariants. We first discuss the many ways one can interpret the dimension vector of a persistence module, with emphasis on those interpretations which are rooted in homological algebra. We then discuss how similar interpretations are used to understand more complicated invariants. We focus in particular on the recent works \cite{KM,AENY,BOO}.

In Section~\ref{sec:relHom}, we give a brief overview of relative homological algebra  following Auslander and Solberg \cite{AS}. We work in this section with an arbitrary finite-dimensional algebra $\Lambda$ and a finite set $\X$ of indecomposable modules. In Section~\ref{sec:approx}, we define the notions of approximations and resolutions by $\X$, and use these to define the \emph{$\X$-dimension} of $\Lambda$ (Definition~\ref{def:xDim}). We then define an ``exact structure'' $\E_\X$ on $\mods\Lambda$, which in particular contains all short exact sequences ending with an approximation by $\X$. Using these short exact sequences, we define the \emph{relative Grothendieck group} $K_0(\Lambda,\X)$ (Definition~\ref{def:grothendieck}).
In Section~\ref{sec:hom_invariants}, we define the \emph{dim-hom} invariants and \emph{homological invariants} relative to $\mathcal{X}$ (Definition~\ref{def:homologicalInvariants}), which are the main objects of study in this paper. Finally, in Section~\ref{sec:projectivisation1}, we use the theory of \emph{projectivization} to give a sufficient condition for when the relative Grothendieck group $K_0(\Lambda,\X)$ is free abelian (Proposition~\ref{cor:grothendieckFree2}). As a consequence, we obtain the following result.

\begin{thm}[Theorem~\ref{thm:hom_dimhom}]\label{thm:mainB}
    Let $\Lambda$ be a finite-dimensional algebra, and let $\mathcal{X}$ be a finite set of indecomposable $\Lambda$-modules which contains the indecomposable projectives. If the algebra
    $\End_\Lambda(\bigoplus_{R \in \mathcal{X}} R)^{op}$
    has finite global dimension, then the dim-hom invariant relative to $\X$ and the homological invariant relative to $\X$ are equivalent.
\end{thm}

In Section~\ref{sec:spreads}, we  study spread modules in more detail. The main result of this section is a concrete description of the Hom-space between spread modules (Proposition~\ref{prop:homSpace}).

In Section~\ref{sec:newInvariants}, we specialize the results of Section~\ref{sec:relHom} to the case of persistence (and more specifically spread) modules. In 
particular, we prove our first main theorem.
\begin{thm}[Theorem~\ref{thm:one_source_finite_gldim}]\label{thm:mainA}
    Let $\poset$ be a finite poset and let $\X$ be a set of connected spread modules over $\incidence$ which contains the indecomposable projectives.
    \begin{enumerate}
        \item If every spread in $\X$ has a unique source, then the association $M \mapsto [M]_\X$ is a homological invariant.
        \item If every spread in $\X$ is supported on an upset of $\poset$, then the association $M \mapsto [M]_\X$ is a homological invariant.
    \end{enumerate}
\end{thm}

 If $\X$ is the set of all spread modules which contain a unique source, then we call the association $[M] \mapsto [M]_\X$ the \emph{single-source homological spread invariant}.  Furthermore, we note that a special case of Theorem~\ref{thm:mainA}(2) can also be deduced from \cite[Theorem~6.12]{miller}. See remark~\ref{rem:miller}.

We then show that allowing $\X$ to contain spread modules with multiple sources sometimes leads to invariants which are not homological (Example~\ref{ex:infinite_gldim}). It remains an open question whether taking $\X$ to be the set of all connected spread modules yields a homological invariant, even in the case where $\poset$ is a 2-dimensional grid.
Finally, in Section~\ref{sec:infinite}, we discuss the existence of algorithms for computing our homological invariants and which aspects of the theory can be extended to infinite posets.

In Section~\ref{sec:otherInvariants}, we compare our new homological invariants to other invariants used in persistence theory. We first show that, if the Hasse quiver of $\poset$ is Dynkin type A, then the barcode is a homological invariant (Proposition~\ref{prop:barCode}). Likewise, we use Theorem~\ref{thm:mainA}, together with results from \cite{BOO}, to show that the rank invariant (Theorem~\ref{thm:boo}) is a homological invariant. We then prove our second main theorem.

\begin{thm}[Theorem~\ref{thm:finer}]\label{thm:mainC}\
Let $\poset$ be a finite poset. Then the single-source homological spread invariant is finer than the rank invariant on $\mods\incidence$. Moreover, these invariants are equivalent if and only if for all $a \in \poset$ the subposet $\{x \in \poset \mid a \leq x\}$ is totally ordered.
\end{thm}

We also explicitly explain why the signed barcode of~\cite{BOO} and the homological invariant relative to the set of interval modules do not conincide, even though both take values in the same Grothendieck group. See Remark~\ref{rem:recVrk}. Finally, we show that the generalized persistence diagram of \cite{KM}, which yields signed approximations of arbitrary persistence modules by spread modules, is not homological with respect to the set of spread modules. We conjecture that more generally it is also not homological relative to any set of indecomposables. See Corollary~\ref{cor:gen_rk_grid} and Conjecture~\ref{conj:not_homological}.

%%%%%%%%%%%%%%%%%%%%%%%%%%%%%%%%%%%%%%%%%%%%%%%%%%%%%%%%
\subsection{Acknowledgements}
This work was supported by NSERC of Canada. A portion of this work was completed while E. H. was affiliated with the Norwegian University of Science and Technology (NTNU). E. H. thanks NTNU for their support and hospitality. The authors are thankful to Rose-Line Baillargeon, Micka\"el Buchet, Barbara Giunti, Denis Langford, Souheila Hassoun, Ezra Miller, Charles Paquette, Job D. Rock, Ralf Schiffler, and Luis Scoccola for their contributions to numerous discussions pertaining to this project. They are also thankful to an anonymous referee for their careful reading of this manuscript and suggestions for improvement.

%%%%%%%%%%%%%%%%%%%%%%%%%%%%%%%%%%%%%%%%%%%%%%%%%%%%%%%%%

\section{Notation and Terminology}\label{sec:terminology}

In this section, we fix notation and terminology for use in the remainder of this paper. We first recall some notation from the introduction. Fix a field $K$. The notation $\dim(-)$ will always mean $\dim_K(-)$. We denote by $\poset$ a poset, which we will assume to be finite unless otherwise stated. We denote by $\incidence$ the ($K$-)incidence algebra of $\poset$ and by $\mods\incidence$ the category of finitely-generated (right) $\incidence$-modules. Unless otherwise stated, the phrase ``an $\incidence$-module" will always refer to an object of $\mods\incidence$. Given $M$ and $N$ two $\incidence$-modules, a ``morphism'' $f:M \rightarrow N$ refers to an $\incidence$-linear map. We denote by $\Hom_{\incidence}(M,N)$, or just $\Hom(M,N)$, the vector space of morphisms from $M$ to $N$. In Section~\ref{sec:relHom}, we will work over an arbitrary finite-dimensional algebra, and thus the symbol $\incidence$ will be replaced with $\Lambda$ while maintaining our conventions.

We identify $\mods\incidence$ with the category of pointwise finite-dimensional (pfd) persistence modules over $\poset$. That is, we consider a module $M \in \mods\incidence$ as a functor from $\poset$ to the category of finite-dimensional $K$-vector spaces. Given $a \leq b \in \poset$, we denote by $M(a)$ and $M(a,b)$ the result of applying $M$ to the object $a$ and to the unique morphism $a\rightarrow b$ in $\poset$, respectively.

We adopt the common practice of identifying modules with their isomorphism classes. In particular, the term ``subcategory'' will always refer to a subcategory which is full and closed under isomorphism. Moreover, any use of the phrase ``the set of ----- modules" could more precisely be replaced with ``the set of isomorphism classes of ----- modules''.

Given a set $\mathcal{X}$ of $\incidence$-modules, we denote by $\add(\mathcal{X})$ the subcategory of $\mods\incidence$ whose objects are (isomorphic to) finite direct sums of the objects in $\mathcal{X}$. That is, a module $M$ is in $\add(\mathcal{X})$ if and only if there is a finite subset $\mathcal{Y} \subseteq \mathcal{X}$ and a tuple of natural numbers $(n_R)_{R \in \mathcal{Y}}$ such that $M \simeq \bigoplus_{R \in \mathcal{Y}} R^{(n_R)}$. Note that by taking $\mathcal{Y} = \emptyset$, we have that $0 \in \add(\X)$. As an example, if $\mathcal{X}$ is the set of indecomposable projective $\incidence$-modules, then $\add(\X)$ is the subcategory consisting of all projective $\incidence$-modules.

\subsection{Spread modules}\label{sec:spreadDefinition} \

In this section, we define a subclass of $\incidence$-modules which we refer to as \emph{spread modules}. As we explain in Remark~\ref{rem:spreadName}, these are often referred to as ``interval modules'' in the persistence theory literature. We will study this class of modules in more details in Section ~\ref{sec:spreads}.

Let $a\in \poset$ and let $B\subseteq \poset$ be a set of incomparable elements. We say $a\leq B$ (respectively $a\geq B$) if there exists $b\in B$ such that $a\leq b$ (respectively $a\geq b$). Sets of incomparable elements of $\poset$ can then be partially ordered by the relation
$$A\leq B \iff \forall_{a \in A}\forall_{b \in B} \ a\leq B \text{ and } A\leq b.$$

Let $Q$ be the Hasse quiver of $\poset$. We call $[a,c]:=\{b \mid a\leq b\leq c\}$ the \emph{interval} from $a$ to $c$. A subset $S$ of $\poset$ is called \textit{connected} if it is a connected part of $Q$; that is, there exist a non-oriented path between any two points of $S$. A subset $S$ of $\poset$ is called \textit{convex} if for any pair $(a,b)\in S$, $[a,b]\subseteq S$.

Clearly, intervals are connected and convex. We generalize the notion of interval with the following definition.

\begin{defn}\label{def:spreads}
Let $A, B \subseteq \poset$ be two sets of incomparable elements with $A \leq B$. The \textit{spread} from $A$ to $B$ is the subset
$$[A, B] := \{ x\; |  \text{ there exists } a\in A, b\in B \text{ such that } a\leq x \leq b\}.$$ We refer to $A$ and $B$ as the sets of \emph{sources} and \emph{targets} of $[A,B]$, respectively.
\end{defn}

If $A = \{a\}$ contains only a single element, we call $[A,B]$ a \emph{single-source spread} and write $[a,B]$ in place of $[\{a\},B]$. \emph{Single-target spreads} are defined and notated analogously. Note that spreads which are both single-source and single-target are precisely intervals (as in the classical poset-theory literature, see e.g. \cite{ec1}). Spreads are also convex. Indeed, if $x,z\in [A,B]$, then there exist $a$ and $b$ such that $a\leq x$ and $z\leq b$, so for any $y\in [x,z]$, we have $a\leq x\leq y\leq z \leq b$, showing it is convex. In fact, all finite convex subsets are spreads, as made precise in the following proposition.

\begin{prop}\label{prop:spreadsConvex}
Let $X$ be a convex subset of $\poset$. Let $A$ be the set of minimal elements of $X$ and $B$ be the set of maximal elements of $X$. Then $X= [A,B]$.
\end{prop}
\begin{proof}
    It is clear that $A$ and $B$ are sets of incomparable elements and that $A \leq B$. Now
    take $y\in X$ to be neither maximal nor minimal. Then there exist elements in $X$ that are respectively bigger and smaller than $y$, say $x_1\leq y \leq z_1$. We then iteratively find bigger and smaller elements, obtaining a sequence
$$x_n\leq ... \leq x_1 \leq y \leq z_1 \leq ... \leq z_m$$
where $x_n$ is minimal and $z_m$ is maximal. This sequence must end in a finite number of steps because $\mathcal{P}$, and thus also $X$, is finite. 
\end{proof}

We also fix the following notation for convenience.

\begin{defn}
    Let $X\subseteq \poset$ be a set of incomparable elements. By Proposition ~\ref{prop:spreadsConvex}, the sets 
    $$[-\infty,X]:=\{y \mid y \leq X\} \qquad [X,\infty]:=\{y \mid X \leq y\}$$
    are spreads. We call them \emph{downsets} and \emph{upsets}, respectively. Note that downsets are sometimes called ``order ideals'', and are characterized by being closed under taking predecessors. Likewise upsets are sometimes called ``order filters'', and are characterized by being closed under taking successors.
\end{defn}

\begin{defn}\label{def:spreadRep}
Given a spread $[A,B]\subseteq \poset$, we call \emph{spread module} the persistence module $M_{[A,B]}$ defined by
\[M_{[A,B]}(p)= 
\begin{cases}
K & \text{ if } p\in [A,B]\\
0 & \text{ if } p\not\in [A,B]
\end{cases}
\quad \quad \quad \quad
M_{[A,B]}(p,q)= \begin{cases}
1_K & \text{ if } p,q \in [A,B]\\
0 & \text{ else}
\end{cases}
\]
If in addition $|A| = 1 = |B|$, we call $M_{[A,B]}$ an \emph{interval module}. We likewise define \emph{single-source spread modules, single-target spread modules, upset modules, and downset modules} in the natural way. 
\end{defn}

\begin{rem}\label{rem:spreadName}
    In much of the recent literature on persistence modules (e.g. \cite{AENY,ABENY,BOO} and several others) the term ``interval module'' is used for what we have called a ``spread module'' and the term ``segment module'' is used for what is called an ``interval module''. Our naming convention was chosen to emphasize that spreads are not intervals in the classical poset-theory sense.
\end{rem}

The following are some important examples of spread modules.

\begin{ex}\label{ex:spreads}\
    \begin{enumerate}
        \item All simple modules are interval modules. Indeed, let $a \in \poset$. Then $S_a := M_{[a,a]}$ is the simple module supported at $a$.
        \item All indecomposable projective modules are upset modules. Indeed, for $a \in \poset$ we have that $P_a:=M_{[a,\infty]}$ is the projective cover of $S_a$. Symmetrically, indecomposable injective modules are downset modules: $I_a:= M_{[-\infty,a]}$ is the injective envelope of $S_a$.
        \item If the Hasse quiver of $\poset$ is Dynkin type $A$, then all indecomposable modules are spread modules. If moreover $\poset$ is totally ordered, then all indecomposable modules are interval modules.
    \end{enumerate}
\end{ex}

\begin{rem}\label{rem:free}
    The modules $P_a$ in Example~\ref{ex:spreads}(2) are sometimes referred to as ``free persistence modules''. This is because if $\poset = \mathbb{Z}^n$, then the module $P_a$ is a free $n$-graded module over $K[x_1,\ldots,x_n]$. See \cite[Section~4.2]{CZ}. On the other hand, $P_a$ is not free as an $\incidence$-module, so we have chosen not to use the word free in this context. See also \cite[Section~6]{BM}.
\end{rem}

By construction, if $M$ is a spread module, then $\dim M(a) \leq 1$ for all $a \in \poset$. Modules with this property are sometimes called \emph{thin}. Another motivation for studying spread modules is the following.

\begin{thm}\cite[Theorem~24]{ABENY}\label{thm:2Dspreads}
Suppose $\poset = [0,a]\times[0,b]\subseteq\mathbb{Z}^2$ is a product of two finite totally ordered sets. Then every thin indecomposable module in $\mods\incidence$ is isomorphic to a spread module.
\end{thm}

As previously mentioned, if $\poset$ is totally ordered, then all indecomposable modules are spread modules. Theorem~\ref{thm:2Dspreads} thus points towards a natural invariant: the multiplicity of thin modules in the direct sum decomposition. If the quiver of $\poset$ is Dynkin type $A$, this invariant is equivalent to the barcode. More generally, however, this invariant loses all information about any indecomposable direct summand which isn't thin. Nevertheless, in some applications, it looks like the proportion of non-thin factors is small, making this invariant somewhat valuable. See e.g. \cite[Section~5]{EH} for further discussion. 

An immediate consequence of Theorem~\ref{thm:2Dspreads} is that, over these particular posets, thin indecomposable modules are uniquely determined by their support. For more general posets, however, there are thin indecomposables which are neither isomorphic to spread modules nor determined by their support, as seen in the following example.

\begin{ex}\label{ex:non_thin}
As in \cite[Example~2.7]{miller}, let $\poset$ be the poset with Hasse diagram
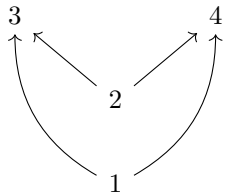
\begin{center}
    \begin{tikzcd}
    3&&4\\
    &2\arrow[ur]\arrow[ul]\\
    &1\arrow[uur,bend right]\arrow[uul,bend left]
    \end{tikzcd}
\end{center}
Note that the incidence algebra $\incidence$ is in this case a path algebra of type $\widetilde{A}_3$. Now for $\lambda \in K$, let $N_\lambda$ be the following thin module.
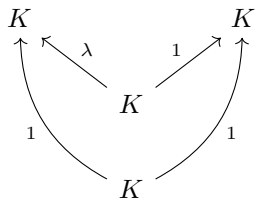
\begin{center}
    \begin{tikzcd}
    K&&K\\
    &K\arrow[ur,"1"above left]\arrow[ul,"\lambda"above right]\\
    &K\arrow[uur,bend right,"1" below right]\arrow[uul,bend left,"1" below left]
    \end{tikzcd}
\end{center}
It is well-known that each $N_\lambda$ is indecomposable and that $N_\lambda \simeq N_{\lambda'}$ if and only if $\lambda = \lambda'$. Moreover, we observe that $N_\lambda$ is a spread module if and only if $\lambda = 1$ (in which case we have $N_\lambda = M_{[\{1,2\},\{3,4\}]}$). This shows that $\mods\incidence$ contains thin indecomposables which are neither isomorphic to spread modules nor determined by their support.
\end{ex}

%%%%%%%%%%%%%%%%%%%%%%%%%%%%%%%

\section{Motivation and related invariants}\label{sec:motivation}

In this section, we examine many invariants from the literature and how the information they contain is interpreted. In particular, we highlight places where homological algebra can be used to clarify, complement, or expand existing frameworks. This section is not meant as a thorough treatise, but rather as motivation for studying the particular invariants we propose.

While we give formal definitions in Definition~\ref{def:invariant}, we recall that an (additive) invariant is a map $p:\mods\incidence \rightarrow \mathbb{Z}^n$ which is constant on isomorphism classes and satisfies $p(M\oplus N) = p(M) + p(N)$. Furthermore, we say that two invariants $p$ and $q$ are \emph{equivalent} if $p(M) = p(N)$ if and only if $q(M) = q(N)$ for all $M, N \in \mods\incidence$.

\subsection{The dimension vector/Hilbert function}\label{sec:dim}
 We begin by overviewing many different ways that one may interpret the dimension vector. Recall that for a persistence module $M$, the \emph{dimension vector} (or \emph{Hilbert function}) is defined by
$$\undim M = (\dim M(a))_{a \in \poset}.$$
This invariant is often considered as taking values in the free abelian group with basis $\mathcal{P}$; however, we can also see $\undim{M}$ as an element of the free abelian group with basis $\{[S_a] \mid a \in \poset\}$. (Precisely, we are treating $\undim M$ as an element of the classical Grothendieck group of $\mods \incidence$. See Section~\ref{sec:relHom}.) The second choice highlights the fact that given a short exact sequence
$$0 \rightarrow L \rightarrow M \rightarrow N \rightarrow 0$$
of persistence modules, one has that
$$\undim(M) = \undim(L) + \undim(N).$$
In particular, the dimension vector $\undim M$ and the multiset of composition factors for $M$ uniquely determine one another. As a special case, we have that $\undim(M\oplus N) = \undim(M) + \undim(N)$; that is, the dimension vector is an (additive) invariant.

We now turn towards understanding the dimension vector $\undim M$ using the projective modules $P_a$. To start with, recall the well known fact that for $a \in \mathcal{P}$, one has $\dim M(a) = \dim\Hom_{\incidence}(P_a,M)$. One may then see $\undim M$ as counting the number of ``test morphisms'' from a set of well-understood modules (namely the indecomposable projectives). Moreover, if $\incidence$ has finite global dimension, then there is a change of basis $\sigma: \mathbb{Z}^{|\poset|} \rightarrow \mathbb{Z}^{|\poset|}$ which can be conceptualized as sending the free abelian group with basis $\{S_a \mid a \in \poset\}$ to that with basis $\{P_a \mid a \in \poset\}$. Formally, $\sigma$ is defined as the inverse of the \emph{Cartan matrix} of $\incidence$. Since we are working over a poset algebra, this also coincides with the M\"obius inversion formula. We will give further details in Section~\ref{sec:relHom}, but for the purpose of this section, the change of basis $\sigma$ works as follows. Given an arbitrary projective module $P$, there is a unique direct sum decomposition $P \simeq \bigoplus_{a \in \poset} (P_a)^{r_a}$.
We then set $[P] = \sum_{a \in \poset} r_a [P_a]$, where $[P_a]$ is the basis element corresponding to $P_a$. Now given $M$ an arbitrary persistence module, we choose a finite projective resolution
$$P_m \rightarrow \cdots \rightarrow P_1 \rightarrow P_0 \rightarrow M.$$
One then has $\sigma(\undim M) = \sum_{j = 1}^m (-1)^j [P_j]$. Taking this one step further, let us denote by $[M]_a$ the coefficient of $[P_a]$ in $\sigma(\undim M)$. Now define
$$P_+ = \bigoplus_{a: [M]_a > 0} (P_a)^{[M]_a}, \qquad\qquad P_- = \bigoplus_{a: [M]_a < 0} (P_a)^{-[M]_a}.$$
One may then consider the pair $(P_+,P_-)$ as a ``signed approximation'' of $M$ by the category of projective modules. Indeed, by applying $\sigma^{-1}$, we obtain the equation $\undim M = \undim P_+ - \undim P_-$.

To summarize, we have the following ways to interpret the dimension vector $\undim M$:
\begin{enumerate}
    \item as recording the dimensions of the vector spaces comprising $M$.
    \item as recording the number of ``test morphisms'' from the indecomposable projectives to $M$.
    \item as some data determined by a projective resolution of $M$.
    \item as a ``signed approximation'' of $M$ by the category of projective modules.
\end{enumerate}

As we will see in Section~\ref{sec:relHom}, ``the homological invariants'' we introduce in this paper are readily given interpretations in the spirit (2), (3), and (4) above. In essence, we will enlarge the set of indecomposable projective modules to a larger set $\mathcal{X}$ of indecomposable $\incidence$-modules, typically a subset of the connected spread modules. This instantly gives rise to the invariant $(\dim\Hom_{\incidence}(R,M))_{R \in \mathcal{X}}$, which can be seen as counting ``test morphisms'' from the objects in $\mathcal{X}$ to $M$. For many choices of $\mathcal{X}$, this is equivalent to data coming from an ``$\mathcal{X}$-resolution'' of $M$. As each step of such a resolution is an ``approximation'' in some precise sense, the result is readily interpreted as a ``signed approximation'' of $M$ by $\mathcal{X}$. Finally, as made explicit in Proposition~\ref{prop:spread_of_persistence}, it is also possible to understand our invariants more directly in the spirit of (1) above.

In the remainder of this section, we focus on many of the invariants which serve as motivation and background for our study. In particular, we emphasize existing works which yeilds interpretations in the spirt of (2), (3), and (4) above. We also discuss alternative approaches to extracting data from a projective resolution, as these can also be readily applied to our new framework.

%%%%%%%%%%%%%%%%%%%%%%%%%%%%%

\subsection{The rank invariant and its generalizations}
Recall the definition of the rank invariant $\unrk M$ from Equation~\ref{eqn:rankInv}. In this section, we explain various interpretations and generalizations of the rank invariant which have appeared in recent work.

\subsection*{The generalized rank invariant and compressed multiplicities}
The rank invariant can be seen as associating one nonnegative integer\footnote{We recall that we are using the term ``interval'' in its classical order-theoretic sense. This is different from the definition used in \cite{KM}.} to $M$ for each interval $I = [a,b] \subseteq \poset$. The recent works \cite[Section~3.2]{thomasThesis}, \cite{KM}, and \cite{AENY} each give alternative interpretations of these integers which extend naturally to allow one to consider subsets of $\poset$ which are not intervals. This results in new invariants, known as the multirank invariant \cite[Section~3.2]{thomasThesis}, the generalized rank invariant \cite{KM}, and compressed multiplicities \cite{AENY}. We give a brief overview of these constructions.

\begin{rem}
    In \cite{KM}, the authors work over general posets (without the assumption that $|\poset| < \infty$), but we keep the assumption that $|\poset| < \infty$ for simplicity. Moreover, they consider functors with target categories more general than $\mathrm{vec}(K)$.
\end{rem}

The definition that follows makes use of the categorical notions of limits and colimits. A brief explanation of these constructions can be found in \cite[Section~3.1 and Appendix~A]{KM}. Note also that our definition is slightly more general than that of \cite{KM}, in which it is required that $X$ be (path-)connected as a subset of $\poset$.

\begin{defn}\label{def:genRank}
    Let $M$ be a ($\poset$-)persistence module and let $X \subseteq \poset$ such that $X$ is a connected poset with the induced order from $\poset$\footnote{Note that $X$ may be a connected poset even if it is not connected as a subset of $\poset$. For example, take $\poset = \{0,1\}\times\{0,1\}$. Then $X = \{(0,0), (1,1)\}$ is connected as a poset, but not as a subset of $\poset$.}. We define a diagram $M(X)$ in the category of finite-dimensional $K$-vector spaces as follows:
    \begin{itemize}
        \item The objects of $M(X)$ are the vector spaces $M(a)$ for $a \in X$.
        \item The morphisms of $M(X)$ are the linear maps $M(a,b)$ for $a \leq b \in X$.
    \end{itemize}
    Alternatively, one may view $M(X)$ as the restriction of $M$ to $X$. We denote by $\rk(M,X)$ the rank of the natural map from the limit of $M(X)$ to the colimit of $M(X)$. We refer to $\rk(M,X)$ as the \emph{$X$-rank} of $M$. More generally, for $\mathcal{R}$ a subset of the power set $2^\poset$ consisting of sets $X$ as in Definition~\ref{def:genRank}, we define the \emph{$\mathcal{R}$-rank} of $M$ to be $\rk(M,\mathcal{R}) := (\rk(M,X))_{X \in \mathcal{R}}$.
\end{defn}

If $X = [a,b]$ is an interval subset of $\poset$, then $\rk(M,X) = \rk(M(a,b))$. In particular, if $\mathcal{I}$ is the set of intervals in $\poset$, then $\rk(M,\mathcal{I}) = \unrk(M)$. Motivated by this fact, Kim and M\'emoli define the \emph{generalized rank invariant} in \cite{KM}. This is precisely the invariant $\rk(M,\mathcal{R})$ for $\mathcal{R}$ the set of ``path connected'' subsets of $\poset$. In the subsequent works \cite{DKM,KM2}, the set $\mathcal{R}$ has been restricted to that of (all) connected spreads. More generally, in \cite{BOO}, Botnan, Oppermann, and Oudot consider the invariants $\rk(M,\mathcal{R})$ where $\mathcal{R}$ is any set of the connected spreads. (They actually work over posets which may not be finite and impose some conditions on these sets, but these conditions are always satisfied over finite posets.)

In another recent work \cite{AENY}, Asashiba, Escolar, Nakashima, and Yoshiwaki work exclusively with finite posets of the form $[0,a] \times [0,b] \subseteq \mathbb{Z}^2$. For every spread $[A,B]$, they introduce three invariants. It is shown in \cite{KM2}, using \cite[Lemma~3.1]{CL}, that the ``total compression factor'' of \cite{AENY} is precisely $\rk(M,[A,B])$. An analogous argument shows that the ``sink-source compression factor'' and ``corner-complete compression factor'' of \cite{AENY} are both of the form $\rk(M,X)$ for some choice of $X$. For example, let $[A,B] \subseteq \poset$ be a connected spread. Then the invariant $d^{ss}(M,[A,B])$ introduced in \cite{AENY} is precisely $\rk(M,A\cup B)$.

Finally, the ``multirank function'' introduced by Thomas in \cite[Section~3.2]{thomasThesis} is similar in concept to the association $X \mapsto \rk(M,X)$, but it does not rely on the formal concepts of limits and colimits. Indeed, given a pair of subsets $X, Y \subseteq \poset$ with no restrictions, Thomas defines the ``multirank from $X$ to $Y$'' as the rank of a map $\oplus_{a \in X} M(a) \rightarrow \oplus_{b \in Y} M(b)$ formed by the component maps of $M$ together with zero maps where necessary. When $X = A$ and $Y = B$ are antichains and $[A,B]$ is a connected spread, we expect that the rank of this map will be related to $\rk(M,A \cup B)$.

%%%%%%%%%%%%%%%%%%

\subsection*{Signed barcodes and M\"obius inversion}
In this section, we will discuss how M\"obius inversion can be used to reinterpret the $\mathcal{R}$-rank invariants. We refer to \cite[Section~3.7]{ec1} for the definition of M\"obius inversion.

As previously mentioned, when $\mathcal{P}$ is totally ordered the (classical) rank invariant is equivalent to the barcode (see also Section~\ref{sec:barcode} for the formal statement). To unpack this fact, let us view the rank invariant and barcode as taking values in the free abelian groups whose bases are given by the sets of (i) interval subsets of $\poset$ and (ii) interval modules over $\incidence$. Then M\"obius inversion is a change of basis $\sigma: \mathbb{Z}^{{|\poset|}\choose{2}} \rightarrow \mathbb{Z}^{{|\poset|}\choose{2}}$ which sends $\unrk(M)$ to the barcode of $M$. We emphasize that this change of basis is \emph{not} the one that sends the interval $[a,b]$ to the interval module $M_{[a,b]}$.

We now return to the case where $\poset$ is not totally ordered. As in the preceding section, let $\mathcal{R}$ be a subset of the power set $2^\poset$ consisting of sets $X$ as in Definition~\ref{def:genRank}. Then $\mathcal{R}$ is itself a poset under inclusion. We denote by $\delta(M,\mathcal{R})$ the result of applying the resulting M\"obius inversion formula to the invariant $\rk(M,\mathcal{R})$. Invariants obtained in this way are sometimes referred to as ``generalized persistence diagrams". See e.g. \cite{patel,KM,MP,BBE,AENY,BE,BOO}.

In \cite{KM}, \cite{BOO}, and \cite{AENY}, the invariants $\delta(M,\mathcal{R})$ are considered for the choices of $\mathcal{R}$ (and generality of poset) discussed in the previous section. To better understand these invariants, given a set $\mathcal{R}$ of connected spreads, we denote $\mathcal{M}(\mathcal{R}) = \{[M_{[A,B]}] | [A,B] \in \mathcal{R}\}$. In all three works, the invariant $\delta(M,\mathcal{R})$ is considered as taking values in the free abelian group which has $\mathcal{M}(\mathcal{R})$ as a basis\footnote{To be precise, we note that the M\"obius inversion procedures for the ``corner-complete'' and ``sink-source'' invariants of \cite{AENY} are done using connected spreads. Since these invariants are both of the form $\rk(M,\mathcal{R})$ for sets $\mathcal{R}$ which contain non-spreads, the M\"obius inversion procedue described here would lead to a different result. Nevertheless, the resulting invariants would still be equivalent since M\"obius inversion is an invertible procedure.}. It is then shown that for any $[A,B] \in \mathcal{R}$ one has $\delta(M_{[A,B]},\mathcal{R}) = [M_{[A,B]}]$. Due to this fact, one can interpret these particular invariants as ``signed approximations'' by spread modules. The word ``signed'' here refers to the fact that some of the coefficients of $\delta(M,\mathcal{R})$ may be negative.

Finally, for $[A,B] \in \mathcal{R}$, let us denote by $r_{[A,B]}$ the coefficient of $[M_{[A,B]}]$ in $\delta(M,\mathcal{R})$. Now define
$$M_+ = \bigoplus_{[A,B]: r_{[A,B]} > 0} (M_{[A,B]})^{r(M,[A,B])}, \qquad\qquad M_- = \bigoplus_{[A,B]: r_{[A,B]} < 0} (M_{[A,B]})^{-r(M,[A,B])}.$$
By applying the inverse M\"obius inversion formula, one then obtains the equation $\rk(M,\mathcal{R}) = \rk(M_+,\mathcal{R}) - \rk(M_-,\mathcal{R})$. See \cite[Theorem~2.5]{BOO}, where the pair $(M_+,M_-)$ is referred to as a ``signed rank decomposition". In the special case that $\mathcal{R} = \mathcal{I}$ is the set of intervals, we recall that $\rk(M,\mathcal{I}) = \unrk(M)$. In this case, the pair $(M_+,M_-)$ admits a visual interpretation called the ``signed barcode''. See \cite[Section~6]{BOO}.

To summarize, there has been a large amount of recent work aimed at turning the $\mathcal{R}$-rank invariants $\rk(M,\mathcal{R})$ into signed approximations $\delta(M,\mathcal{R})$ by spread modules. One of our goals is to examine if and when these coincide with the homological notion of a (signed) approximation by spread modules, as defined in Section~\ref{sec:relHom}. For the (classical) rank invariant, \cite{BOO} has already given an affirmative answer to this question, which we will now explain.

%%%%%%%%%%%%%%%%%%%%%%%%%%%%%%%%
\subsection*{Test morphisms and exact structures}
In \cite{BOO}, the authors consider a set $\mathcal{H}$ of spreads called \emph{hooks}. (See Section~\ref{sec:rank} for the definition.) They then show that the rank invariant $\unrk(M)$ is equivalent to the (dim-hom) invariant $(\dim_K\Hom_\incidence(H,M))_{H \in \mathcal{H}}$. In other words, the rank invariant can be seen as recording the number of ``test morphisms'' from the hook modules to $M$.

It is further shown that the hook modules are precisely the indecomposable ``relative projective'' modules for some \emph{exact structure} on $\mods\incidence$. We will discuss exact structures in more detail in Section~\ref{sec:relHom}. For the purposes of motivation, this amounts to choosing a subset of short exact sequences in $\mods\incidence$ which are deemed ``admissible''. In \cite{BOO}, these are the short exact sequences $0 \rightarrow L \rightarrow M \rightarrow N \rightarrow 0$ for which $\unrk(M) = \unrk(L) + \unrk(N)$. As a consequence of our Theorem~\ref{thm:finer}, this means that the rank invariant is equivalent to a ``signed approximation'' by hook modules in the homological algebra sense. In other words, the rank invariant is a ``homological invariant'' as defined in Definition~\ref{def:homologicalInvariants}. In Section~\ref{sec:rank}, we more carefully examine how this compares to the signed approximation coming from the invariant $\delta(M,\mathcal{H})$.

%%%%%%%%%%%%%%%%%%%%%%%%%%%%%%

\subsection{Generalized Betti numbers and $g$-vectors}
As the dimension vector can be derived from the data of a projective resolution, so too can homological invariants be derived from the data of a ``relative projective resolution". In our discussion of the dimension vector, this is done by taking the signed sum of the terms which appear in the resolution. However, this is not the only way to turn the data of a resolution into a vector which is readily usable in machine learning. We will briefly explain some of the alternative vectorization methods, all of which can be used in a straightfoward manner for our relative homological algebra framework.

Let $M$ be an $\incidence$-module and assume it admits a finite projective resolution
$$P_m \xrightarrow{f_m} \cdots \xrightarrow{f_2} P_1 \xrightarrow{f_1} P_0 \xrightarrow{f_0} M.$$
We will further assume that this projective resolution is taken to be ``minimal'', meaning that each map $f_k$ is the projective cover of $\ker(f_{k-1})$. For $a \in \poset$, and $k \in \{0,\ldots,m\}$, we denote by $r_{k,a}$ the multiplicity of $P_a$ as a direct summand of $P_k$. Each $r_{k,a}$ can the be seen as an invariant of $M$, and in particular each tuple $(r_{k,a})_{a \in \poset}$ is referred to as the \emph{$k$-th generalized Betti numbers} of $M$. This is a straightforward analog of classical Betti numbers, where one records the multiplicities of free modules in a free resolution.

We recall that the dimension vector can be recovered from the generalized Betti numbers. Indeed, the dimension vector $\undim M$ is then (equivalent to) $\sum_{k = 0}^m\sum_{a \in \poset}r_{k,a} [P_a]$, where $[P_a]$ is the basis element associated to $P_a$. This means that the generalized Betti numbers contain more information than the dimension vector. On the other hand, we caution that, unlike the dimension vector, the generalized Betti numbers require the assumption that our chosen projective resolution is minimal. Without this, they will no longer be well defined.

Especially in cases where an algebra has infinite global dimension, one may also record only the first $k$ generalized Betti numbers for some choice of $k$. This is related to the concept of ``$g$-vectors'' from representation theory. Given a module $M$, its \emph{$g$-vector} is by definition $(r_{0,a}-r_{1,a})_{a \in \poset}$. This invariant is particularly important to the concepts of cluster algebras \cite{FZ4}, $\tau$-tilting theory \cite{AIR}, and stability conditions \cite{king}.

%%%%%%%%%%%%%%%%%%%%%%%%

\subsection{Amplitudes and a view towards stability}

We conclude our motivation section with a brief discussion of stability. In this paper, we primarily discuss invariants as a tool for determining whether two persistence modules belong to the same isomorphism class. However, in applications, one is often interested in whether two persistence modules are derived from similar data (or from similar filtrations of topological spaces). In order for invariants to be useful for this purpose, one turns to \emph{stability}. Given an invariant $p$, the key idea is to define a meaningful notion of distance on the output space of $p$. Informally, one then says that $p$ is stable if given two ``similar'' persistence modules $M$ and $N$, one has that the distance between $p(M)$ and $p(N)$ is sufficiently small.

After the first version of this manuscript was posted to arXiv, Oudot and Scoccola proved a stability result for certain invariants coming from exact structures \cite[Theorem~26]{OS}. Their result is proven with respect to the ``Bottleneck dissimilarity function" for finitely presented $\mathbb{R}^n$-persistence modules. Another place one may turn to look for stability results is the recent work of Giunti, Nolan, Otter, and Waas on \emph{amplitudes} and their stability \cite{amplitudes}. Amplitudes are $\mathbb{Z}$-valued invariants which are defined for any abelian category, and are characterized by how they behave with respect to short exact sequences. One could readily adapt this definition to the generality of exact categories by considering only those short exact sequences which are deemed ``admissible''. It remains an open question to determine whether such a modification could be used to prove stability results for the invariants introduced in this paper.

%%%%%%%%%%%%%%%%%%%%%%%%%%%%%%%%%%%%%%%%%%%%%%%%%%%%%%%%

\section{Relative homological algebra}\label{sec:relHom}

The purpose of this section is to provide background on relative homological algebra and relative Grothendieck groups. We largely follow the setup of Auslander and Solberg \cite{AS}. As these results are fully general, throughout this section we will work over a finite-dimensional $K$-algebra $\Lambda$ and assume only (i) that $\Lambda$ is basic and (ii) that every simple $\Lambda$-module is 1-dimensional over $K$. Nevertheless, we note that the incidence algebras $\incidence$ satisfy both (i) and (ii), and so readers interested mainly in applications to persistence theory can safely replace $\Lambda$ with $\incidence$.

The assumptions (i) and (ii) are made to facilitate the discussion of the dimension vector in the sequel. We recall its more general definition now.

\begin{defn}
    Let $\mathcal{U}$ be the set of indecomposable projective $\Lambda$-modules, and let $M \in \mods\Lambda$. Then the \emph{dimension vector} of $M$ is
    $$\undim M := (\dim_K\Hom_\Lambda(P,M))_{P \in \mathcal{U}}.$$
\end{defn}

Alternatively, the assumptions (i) and (ii) are equivalent to the existence of a finite quiver $Q$ such that $\Lambda$ is a quotient of the path algebra $KQ$ by an admissible ideal. In particular, every module $M \in \mods\Lambda$ can be considered as a representation of $Q$ analogously to the way an $\incidence$-module is seen as a functor from $\poset$. For each vertex $v$ of $Q$, we denote by $M(v)$ the vector space at vertex $v$ arising from viewing $M$ in this way. There is then a bijection $\sigma$ from $\mathcal{U}$ to the set of vertices of $Q$ such that for all $P \in \mathcal{U}$ one has $\dim_K M(\sigma(P)) = \dim_K\Hom_{\Lambda}(P,M)$.

%%%%%%%%%%%%%%%%%%%%%%%%%%%%%%%%%%%%
\subsection{Approximations, exact structures, and Grothendieck groups}\label{sec:approx}

Let $\mathcal{X}$ be a finite set of indecomposable modules, and let $M$ be an arbitrary module. We recall that $\add(\X)$ is the subcategory of $\mods\Lambda$ consisting of all finite direct sums (including the empty direct sum) of modules in $\X$. We say a morphism $f: R\rightarrow M$, where $R \in \add(\X)$, is a (right) \emph{$\X$-approximation} if every morphism $R' \rightarrow M$ with $R' \in \add(\X)$ factors through $f$. Since we have assumed $\X$ to be finite, such a morphism is guaranteed to exist. We call $f$ a \emph{minimal} approximation by $\X$ if for all $g: R \rightarrow R$ such that $f \circ g = f$, the endomorphism $g$ is an isomorphism. This is equivalent for asking the domain of $f$ to have the minimum possible number of direct summands.

The following lemma is crucial to the remainder of our setup.

\begin{lem}\label{lem:epiApprox}
    Suppose $\X$ contains all of the indecomposable projective modules and let $f: R \rightarrow M$ be an approximation by $\X$. Then $f$ is an epimorphism.
\end{lem}

\begin{proof}
     Let $q: P \rightarrow M$ be the projective cover of $M$. By assumption, $P \in \add(\X)$, so there exists $g: P \rightarrow R$ such that $q = f \circ g$. Since $q$ is an epimorphism, this implies that $f$ is an epimorphism as well.
\end{proof}

We now consider two examples. Note in particular that Example~\ref{ex:projCoverBad}(2) below shows that there may exist epimorphisms from $\add(\X)$ to $M$ which are not $\X$-approximations.

\begin{ex}\label{ex:projCoverBad} \
    \begin{enumerate}
        \item If $\X$ consists only of the indecomposable projectives, then any epimorphism $f: R \rightarrow M$ with $R \in \add(\X)$ is an $\X$-approximation. Indeed, any morphism $g: R'\rightarrow M$ with $R' \in \add(\X)$ will factor through $f$ by the definition of projective. Moreover, such an epimorphism $f$ will be minimal if and only if it is the projective cover of $M$.
        \item Let $\Lambda = K(1\rightarrow 2)$. Then $\mods\Lambda$ contains three indecomposable modules: the projective-injective $P_1$, the simple-projective $P_2$, and the simple-injective $S_1$. Let $\X = \{P_1,P_2,S_1\}$. Then the identity $S_1\rightarrow S_1$ is an $\X$-approximation which does not factor through the projective cover $P_1 \rightarrow S_1$. This shows that the projective cover $P_1\to S_1$, which is an epimorphism from add$(\X)$, is not an $\X$-approximation.
    \end{enumerate}
\end{ex}

We will assume from now on that $\X$ contains all of the indecomposable projectives. In light of the previous lemma, we can then consider an analogue of a projective resolution of $M$, formed by taking approximations by $\X$ at each step. This yields the following diagram, where each $p_j$ is an approximation by $\X$ and $q_j=i_{j-1} \circ p_j$ for every $j$.

\begin{center}
    \begin{tikzcd}[column sep = 4em,row sep = 3em]
        \cdots \arrow[r]\arrow[dr] & R_2 \arrow[r,"q_2"]\arrow[dr,"p_2"] & R_1 \arrow[r,"q_1"]\arrow[dr,"p_1"] & R_0 \arrow[r,"q_0"]\arrow[dr,"p_0"] & M \arrow[r,"q_{-1} \ = \ 0"] & 0\\
        &\ker(q_2)\arrow[u,"i_2"right]&\ker(q_1)\arrow[u,"i_1"right]&\ker(q_0)\arrow[u,"i_0"right]&\ker(q_{-1}) = M \arrow[u,"i_{-1} = 1_M"right]
    \end{tikzcd}
\end{center}
The first row of this diagram is exact by definition, and is called a \emph{resolution of $M$ by $\X$} or \emph{$\X$-resolution}. If every $p_i$ is a minimal approximation, then the resolution itself is called minimal. We then have the following definitions (see \cite[Section 8.4]{enochsJenda}).

\begin{defn}\label{def:xDim}\
    \begin{enumerate}
        \item Let $M \in \mods\Lambda$ and let $R_\bullet$ be a minimal resolution of $M$ by $\X$. The \emph{$\X$-dimension} of $M$ is then defined as $\min(\{j \mid R_{j+1} = 0\} \cup \{\infty\})$.
        \item The \emph{$\X$-dimension} of $\Lambda$ is defined to be the supremum over $M \in \mods\Lambda$ of the $\X$-dimension of~$M$.
        \item If there exists a module of infinite $\X$-dimension, we say the $\X$-dimension of $\Lambda$ is \emph{properly infinite}. On the other hand, if the $\X$-dimension of $\Lambda$ is infinity but every $M \in \mods\Lambda$ has finite $\X$-dimension, we say that the $\X$-dimension of $\Lambda$ is \emph{effectively infinite}.
    \end{enumerate}
\end{defn}

\begin{rem}\
    \begin{enumerate}
        \item If $\X$ contains only the indecomposable projectives, the $\X$-dimensions of a module $M$ and of the algebra $\Lambda$ coincide with the projective dimension and global dimension, respectively. In this case, it is well-known that algebras with effectively infinite $\X$-dimension do not exist.
        \item The relationship between the projective dimension and $\X$-dimension of a module, and by extension between the global dimension and $\X$-dimension of an algebra, is not clear. For example, in the setup of Example~\ref{ex:projCoverBad}, the projective dimension of $S_1$ is 1, while its $\X$-dimension is 0. On the other hand, in Example~\ref{ex:infinite_gldim} we give an example of a module with projective dimension 1 and infinite $\X$-dimension. One way to explain this discrepancy is that as we increase the size of $\X$, we have both added more modules with which to form $\X$-approximations, but also added more maps which need to factor through these approximations.
    \end{enumerate}
\end{rem}

Now let $\E_\mathrm{max}$ be the class of short exact sequences in $\mods\Lambda$. We then denote
$$\E_\X = \left\{\left(0 \rightarrow M \xrightarrow{f} N \xrightarrow{g} L \rightarrow 0\right)\in \E_{max} \ \middle| \ (\forall R \in \X) \ \Hom(R,g) \text{ is surjective}\right\}.$$
Equivalently, the short exact sequences in $\E_\X$ are precisely those on which the functor $\Hom(R,-)$ is exact for all $R \in \X$. 

It is shown in \cite{DRSS} that $\E_\X$ gives an \emph{exact structure} on $\mods\Lambda$. Readers are referred to \cite{B11} for more information on exact structures. In the present paper we will only use the following straightforward facts about~$\E_\X$.

\newpage

\begin{prop}\label{prop:exactStructure} \
    \begin{enumerate}
        \item Let $0 \rightarrow M \rightarrow N \rightarrow L \rightarrow 0$ be an exact sequence in $\E_\X$. Then
        $$\dim \Hom(R, N) = \dim\Hom(R, L) + \dim\Hom(R, M)$$
        for all $R \in \add(\X)$.
        \item Let $f: R \rightarrow M$ be an approximation of some module $M$ by $\X$. Then
        $$0 \rightarrow \ker(f) \rightarrow R \xrightarrow{f} M \rightarrow 0$$
        is in $\E_\X$. Moreover, $\E_\X$ is the smallest exact structure which contains all such short exact sequences.
        \item Every split short exact sequence is included in $\E_\X$.
    \end{enumerate}
\end{prop}

Bearing this result in mind, we are prepared for the main definition of this section.

\begin{defn}\label{def:grothendieck}\
    \begin{enumerate}
        \item Let $F$ be the free abelian group generated by the symbols $[M]$ for every isomorphism class of $M \in \mods\Lambda$. Denote by $H_\X$ the subgroup generated by the collection
$$\{\; [M]-[N]+[L]\; | \text{ there exists a short exact sequence } 0\to M\to N\to L\to 0 \text{ in }\mathcal{E}_\X\}.$$
    The quotient group $K_0(\Lambda,\X):= F/H_\X$ is called the \emph{Grothendieck group} of $\Lambda$ relative to $\X$.
        \item Given $M \in \mods\Lambda$, we denote by $[M]_\X$ (or just $[M]$ if the set $\X$ is clear from context) the representative of $M$ in the Grothendieck group $K_0(\Lambda,\X)$.
    \end{enumerate}
\end{defn}

Before proceeding, we highlight two important examples of relative Grothendieck groups.

\begin{ex}\label{ex:split}\
    \begin{enumerate}
        \item Let $\X$ be the set of indecomposable projective modules. Then $\E_\X = \E_{max}$ contains all short exact sequences. In this case, we call $K_0(\Lambda) := K_0(\Lambda,\X)$ the \emph{classical Grothendieck group}. This group has a basis given by the simple modules. By our assumptions (i) and (ii) from the beginning of the Section~\ref{sec:relHom}, we can identify $[M] \in K_0(\Lambda)$ with the dimension vector of $M$ using this basis. If in addition we assume that the global dimension of $\Lambda$ is finite, then $K_0(\Lambda)$ also has a dual basis given by $\X$.
        \item Let $\X$ be the set of all indecomposable $\Lambda$-modules. While this will generally not be a finite set, much of the theory developed in this section still applies. In particular, the exact structure $\E_\X$ contains only the split short exact sequences. This exact structure is often denoted $\E_{min}$, and the corresponding relative Grothendieck group is often denoted $K_0^{split}(\Lambda):= K_0(\Lambda, \E_{min})$. This is the free abelian group with basis the set of isomorphism classes of indecomposable $\Lambda$-modules (see e.g. Proposition~\ref{prop:GrothendieckFree} below). That is, given a module $M \in \mods\Lambda$, the data $[M] \in K_0^{split}(\Lambda)$ is equivalent to the direct sum decomposition of $M$. In particular, the isomorphism class of $M$ is uniquely determined by $[M]$.
    \end{enumerate}
\end{ex}

\begin{prop}\label{prop:GrothendieckFree}
If the $\X$-dimension of $\Lambda$ is not properly infinite, the Grothendieck group $K_0(\Lambda,\X)$ is free abelian with basis $\{[R]\mid R\in\X \}$. In particular, $K_0(\Lambda, \X) \simeq \mathbb{Z}^{|\X|}$. Moreover, if $$0\to R_n\xto{q_n}\cdots \xto{q_2} R_1\xto{q_1} R_0\xto{q_0} M\to 0$$ is a finite $\X$-resolution of $M \in \mods\Lambda$, then $[M] = \sum_{i = 0}^n (-1)^i [R_i]$.
\end{prop}

\begin{proof}
Let $M\in\mods\Lambda$ and choose a finite $\X$-resolution of $M$ as in the statement. We can extract short exact sequences
$$0\to \text{ker}(q_i) \into R_i \onto \text{ker}(q_{i+1})\to 0.$$
Each of these short exact sequences is in $\E_\X$ by Proposition~\ref{prop:exactStructure}.
So, in the Grothendieck group, $[R_i]=  [\text{ker}(q_i)] + [\text{ker}(q_{i+1})]$ for all $i$. This yields

\begin{align*}
[M]&= [R_0] - [\text{ker}(q_0)]\\
&= [R_0] - \big( [R_1] -  [\text{ker}(q_1)]\big) \\
&= [R_0] - [R_1] + \big( [R_2] - [\text{ker}(q_2)]\big)\\
&=\cdots\\
&=  \sum_{i=0}^n (-1)^i [R_i],
\end{align*}
where each $R_i$ is in $\text{add}(\X)$.  This shows $\{[R] \mid R\in\X \}$ generates the Grothendieck group. 

It remains to show that $K_0(\Lambda,\X) = \langle [R]\rangle_{R \in \X}$ is free abelian.
To see this, we recall that the short exact sequences ending in an $\X$-approximation generate the subgroup $H_\X$ defining $K_0(\Lambda,\X)$. See Definition~\ref{def:grothendieck} and \cite{AS}. Moreover, given $R \in \add(\X)$, the identity $R\rightarrow R$ is a minimal $\X$-approximation. Together, these facts imply that every exact sequence generating $H_\X$ on the subcategory $\add(\X)$ is split. This concludes the proof.
\end{proof}

%%%%%%%%%%%%%%%%%%%%%%%%%%%%%
\subsection{Homological and dim-hom invariants}\label{sec:hom_invariants}

In this section, we define \emph{dim-hom invariants} and \emph{homological invariants}. In order to make these definitions precise, we first formalize what we mean by an \emph{invariant}, as outlined in the introduction.

\begin{defn}\label{def:invariant}\
    \begin{enumerate}
        \item An \emph{invariant} is a surjective group homomorphism $p: K_0^{split}(\Lambda) \rightarrow \mathbb{Z}^n$ for some $n \in \mathbb{N}$. We call $n$ the \emph{rank} of the invariant $p$.
        \item Given $p$ and $q$ two invariants, we say that $p$ is \emph{finer} than $q$ if $\ker p \subseteq \ker q$. Likewise, we say that $p$ is \emph{equivalent} to $q$ if $\ker p = \ker q$. In particular, $p$ is finer than $q$ if and only if $q(M)=q(N)$ implies $p(M)=p(N)$.
    \end{enumerate}
\end{defn}

To study different types of invariants we consider inclusions of exact structures, see \cite{BHLR}. We need in particular the following lemma.

\begin{lem}\label{lem:finer}
    Let $\X$ and $\mathcal{Y}$ be finite sets of indecomposable modules which contain the indecomposable projectives, and suppose that $\X \subseteq \mathcal{Y}$. Then there is a well-defined quotient map $p: K_0(\Lambda,\mathcal{Y}) \rightarrow K_0(\Lambda,\X)$.
\end{lem}

\begin{proof}
    It is immediate from the definitions of $\E_\X$ and $\E_\mathcal{Y}$ that $\E_\X \supseteq \E_\mathcal{Y}$. The result then follows immediately from the definitions of $K_0(\Lambda,\X)$ and $K_0(\Lambda,\mathcal{Y})$.
\end{proof}

    We note that, even if $\mods\Lambda$ contains infinitely many indecomposable objects, the statement and proof of Lemma~\ref{lem:finer} can be adapted to allow one to take $\X$ to be the set of all indecomposable modules. Thus for any finite set $\mathcal{Y}$ of indecomposable modules containing the indecomposable projectives, we have that $K_0(\Lambda, \mathcal{Y})$ is a quotient of $K_0^{split}(\Lambda)$. This leads to the following definition.
   
\begin{defn}\label{def:homologicalInvariants}
    Let $p: K_0^{split}(\Lambda)\rightarrow \mathbb{Z}^n$ be an invariant, and let $\X$ be a finite set of indecomposable modules.
    \begin{enumerate}
        \item We say that $p$ is a \emph{dim-hom invariant relative to $\X$} if it is equivalent to the invariant given by the map 
        \begin{align*}
        K_0^{split}(\Lambda) &\rightarrow \mathbb{Z}^{|\X|}    \\
        M &\mapsto (\dim\Hom_{\Lambda}(R,M))_{R \in \X}
        \end{align*}
        
        \item     We say that $p$ is \emph{homological relative to} $\X$ if all of the following hold:
    \begin{enumerate}
        \item All of the indecomposable projectives are in $\X$.
        \item The $\X$-dimension of $\Lambda$ is not properly infinite.
        \item $p$ is equivalent to the invariant given by the quotient map $K_0^{split}(\Lambda) \rightarrow K_0(\Lambda, \X)\simeq \mathbb{Z}^{|\X|}.$

    \end{enumerate}
        \item We say that $p$ is a \emph{dim-hom invariant} (resp. \emph{homological invariant}) if it is a dim-hom (resp. homological) invariant relative to some $\X$.
    \end{enumerate}
\end{defn}

 We will show that if a given $\X$ respects some assumptions extending (2b) above, then the homological and dim-hom invariants relative to $\X$ coincide. See Theorem ~\ref{thm:hom_dimhom}.

The following example is well-known.

\begin{prop}\label{prop:dimEquiv}
       Suppose that $\Lambda$ has finite global dimension and let $\X$ be the set of indecomposable projective modules. Then the dimension vector $M \mapsto \undim M := (\dim\Hom_\Lambda(R,M))_{R \in \X}$ and the homological invariant $p: K_0^{split}(\Lambda) \rightarrow K_0(\Lambda)$ are equivalent.
\end{prop}

\begin{rem} In Proposition~\ref{prop:dimEquiv}, one often interprets $\undim M$ as sitting inside of $K_0(\Lambda)$ by taking the simple modules as a basis. One can change to the basis of $K_0(\Lambda)$ given by the indecomposable projectives by sending each simple module to the alternating sum of the terms in its projective resolution. See Section~\ref{sec:dim} for additional discussion.
\end{rem}

We will explore homological invariants (and dim-hom invariants) in the context of persistence modules in Section~\ref{sec:newInvariants}.

%%%%%%%%%%%%%%%%%%%%%%%%%%%
\subsection{Projectivisation}\label{sec:projectivisation1}

Fix a finite set $\mathcal{X}$ of indecomposable modules which contains the indecomposable projectives. In this section, we use the theory of \emph{projectivisation} to construct an algebra with indecomposable projective modules corresponding to modules in $\X$. In particular, this will allow us to show that, for many choices of $\X$, the homological and dim-hom invariants relative to $\X$ coincide. It will also allow us to reduce the problem of computing $\X$-resolutions to that of computing projective resolutions over an endomorphism algebra.
The ideas outlined here are well-known, and readers are referred to \cite[Section~II.2]{ARS} for details.

We recall that a morphism $f: M \rightarrow N$ is called a \emph{section} (resp. \emph{retraction}) if there exists a morphism $g: N \rightarrow M$ such that $g\circ f$ (resp. $f \circ g$) is the identity.

\begin{defn}\label{def:irreducible}
A morphism $f:R \rightarrow R'$, with $R, R'\in \X$, is called $\X$\textit{-irreducible} if it is not itself a section or a retraction and if all factorizations of $f$ of the form $R\xto{g} M\xto{h} R'$, with $M\in \add(\X)$, must have $g$ be a section or $h$ be a retraction.
\end{defn}

Now consider the module $T = T_\X: = \bigoplus_{R\in \X} R$. Then $\Gamma = \Gamma_\X := \End_{\Lambda}(T)^{op}$ is a finite-dimensional $K$-algebra under the composition of morphisms and the standard vector space structure. Denote by $Q$ the Gabriel quiver of $\Gamma$. The vertices of $Q$ can be identified with $\X$. Moreover, for $R, R' \in \X$ the number of arrows $R \rightarrow R'$ in $Q$ coincides with the number of linearly independent $\X$-irreducible morphisms $R' \rightarrow R$ in $\Hom_{\Lambda}(X,Y)$.

The association $M\mapsto \Hom_{\Lambda}(T,M)$ yields a functor from $\mods\Lambda$ to $\mods \Gamma$. The action of $\Gamma$ on $\Hom_{\Lambda}(T,M)$ is given by precomposition. This functor induces an additive equivalence between $\add(\X)$ and the category of projective $\Gamma$-modules. That is, given an (indecomposable) module $R\in \X$, the $\Gamma$-module $\Hom_{\Lambda}(T,R)$ is isomorphic to the projective module (over $\Gamma$) at the vertex corresponding to $R$ in the quiver $Q$.

The following is critical.

\begin{prop}\cite[Proposition~2.1]{ARS}\label{prop:homSame}
    Let $R \in \add(\X)$ and $M \in \mods\Lambda$. Then the functor $\Hom_\Lambda(T,-)$ induces an isomorphism
    $$\Hom_\Lambda(R,M) \simeq \Hom_\Gamma(\Hom_\Lambda(T,R),\Hom_\Lambda(T,M)).$$
\end{prop}

Putting this together, we have the following.

\begin{prop}\label{prop:finiteGlDim}\
\begin{enumerate}
    \item Let $M \in \mods\Lambda$ and let $R_\bullet$ be a (minimal) resolution of $M$ by $\X$. Then $\Hom_{\Lambda}(T,R_\bullet)$ is a (minimal) projective resolution of the $\Gamma$-module $\Hom_{\Lambda}(T, M)$.
    \item Let $M \in \mods\Lambda$ and let $Q_\bullet$ be a minimal projective resolution of the $\Gamma$-module $\Hom_\Lambda(T,M)$. Then there exists a resolution $R_\bullet$ of $M$ by $\X$ such that $\Hom_\Lambda(T,R_\bullet) \simeq Q_\bullet$.
    \item The $\X$-dimension of $\Lambda$ is bounded by the global dimension of  $\Gamma$. In particular, if $\Gamma$ has finite global dimension, then $\Lambda$ has finite $\X$-dimension.
\end{enumerate}
\end{prop}

\begin{proof}
     (1) and (2) follow immediately from Proposition~\ref{prop:homSame} and the fact that $\Hom_{\Lambda}(T,-)$ induces an additive equivalence. (3) is an immediate consequence of (1).
\end{proof}

While the global dimension of $\Gamma$ gives an upper bound on the $\X$-dimension, we demonstrate with the following example that this bound will in general not be tight.

\begin{ex}\label{ex:glDimNotTight}
 Let $\Lambda$ be as in Example~\ref{ex:projCoverBad}. Then $\Gamma = \End_\Lambda(P_1\oplus S_1\oplus P_2)^{op}$ is the so-called \emph{Auslander algebra} of $\Lambda$. The Gabriel quiver of $\Gamma$ is
    $$\begin{tikzcd} Q = & P_2 \arrow[r,"\iota"]& P_1 \arrow[r,"\pi"] & S_1,\end{tikzcd}$$
    and we can write $\Gamma \simeq KQ/(\iota\pi)$. It is well-known that the global dimension of $\Gamma$ is 2, while the $\X$-dimension of $\Lambda$ is clearly 0  since $\add(\X) = \mods\Lambda$.
\end{ex}

As a consequence of Propositions~\ref{prop:finiteGlDim} and~\ref{prop:GrothendieckFree}, we have the following.

\begin{prop}\label{cor:grothendieckFree2}
Suppose that $\Gamma$ is of finite global dimension. Then
\begin{enumerate}
    \item The Grothendieck group relative to $\X$ is free with basis given by $\{[R]\mid R\in\X \}$. In particular, the quotient map $K_0^{split}(\Lambda) \rightarrow K_0(\Lambda,\X)$ is a homological invariant.
    \item There is an isomorphism of Grothendieck groups $K_0(\Lambda,\X) \rightarrow K_0(\Gamma)$ which sends $[M]$ to $[\Hom_\Lambda(T,M)]$.
\end{enumerate}
\end{prop}

\begin{rem}\label{rem:exactGroth} \
    
    \begin{enumerate}
        \item The axioms of exact structures allow one to define a relative version $\mathcal{D}^b( \Lambda, \X)$ of the derived category of $\Lambda$ relative to $\X$, see \cite[Chapter~4]{krause}. 
    When the $\X-$dimension of $\Lambda$ is finite, then it can be shown that $T$ is a tilting object as defined in \cite[Section~7.2]{krause}. In this case the functor $\Hom_{\Lambda}(T,-)$ induces an additive equivalence $\mathcal{D}^b( \Lambda, \X) \to \mathcal{D}^b( \Gamma),$ and thus also an isomorphism of the corresponding Grothendieck groups $K_0(\Lambda,\X) \simeq K_0(\Gamma).$
    It is well-known that the Grothendieck group $K_0(\Gamma)$ is freely generated by the projective $\Gamma$-modules since $\Gamma$ has finite global dimension. 
    As the projectives correspond to the elements of $\X$ under the above isomorphism, we can deduce Proposition~\ref{prop:GrothendieckFree} from the general theory of relative derived categories.
    
    Note however that  the simple $\Gamma$-modules in general do not correspond to $\Lambda$-modules, but only to objects in the relative derived category. 
    So even if the derived categories are equivalent, the $\X$-dimension of $\Lambda$ is in general not equal to the global dimension of $\Gamma$, like it is often observed in classical tilting theory.
        \item We are not aware of any case where the $\X$-dimension of $\Lambda$ is finite and the global dimension of $\Gamma$ is infinity. It would not be surprising if such a case were impossible in light of the fact that tilting preserves finite global dimension, see \cite{KK}.
    \end{enumerate}
\end{rem}

We are now prepared to prove the main results of this section.

\begin{prop}\label{prop:preHom}
    Let $\mathcal{X}$ be a finite set of indecomposable $\Lambda$-modules which contains the indecomposable projectives, and suppose that the global dimension of $\Gamma = \Gamma_\X$ is finite. Given $M \in \mods\Lambda$, write $[M] \in K_0(\Lambda,\X)$ in the form
        $$[M] = \sum_{R \in \X} c_R[R],$$
    where each $c_R \in \mathbb{Z}$. Then for $R' \in \X$, we have
    $$\dim\Hom_\Lambda(R',M) = \sum_{R \in \X}c_R\dim\Hom_\Lambda(R',R).$$
\end{prop}

\begin{proof}
Let $n$ be the $\X$-dimension of $M$. We proceed by induction of $n$.
     
     If $n = 0$, then $M \in \add(\X)$. The result then follows immediately from the additivity of the functors $\Hom(R',-)$.
     
     Now let $n > 0$ and suppose the result holds for $n - 1$. Let $f: R \twoheadrightarrow M$ be a minimal $ \X$-approximation of $M$. Note that $f$ is surjective by Lemma~\ref{lem:epiApprox} and the short exact sequence
     $$0 \rightarrow \ker f \rightarrow R \rightarrow M \rightarrow 0$$
     is in $\mathcal{E}_\X$ by Proposition~\ref{prop:exactStructure}(2). So we have $[M] = [R] - [\ker f]$. Moreover, for $R' \in \X$, Proposition~\ref{prop:exactStructure}(1) implies that
     $$\dim\Hom_\Lambda(R', M) = \dim\Hom_\Lambda(R', R) - \dim\Hom_\Lambda(R',\ker f).$$
     Since the $\X$-dimensions of $R$ and $\ker f$ are 0 and $n-1$, respectively, the result then follows from the induction hypothesis.
\end{proof}

\begin{thm}[Theorem~\ref{thm:mainB}]\label{thm:hom_dimhom}
  Let $\mathcal{X}$ be a finite set of indecomposable $\Lambda$-modules which contains the indecomposable projectives, and suppose that the global dimension of $\Gamma = \Gamma_\X$ is finite.
    Let $p: K_0^{split}(\Lambda) \rightarrow K_0(\Lambda,\X)$ be the canonical quotient map, and let $q: K_0^{split}(\Lambda) \rightarrow \mathbb{Z}^{|\X|}$ be the dim-hom invariant given by $q(M) = (\dim\Hom_\Lambda(R,M))_{R \in \X}$. Then $p$ and $q$ are equivalent invariants. In particular, any invariant that is either homological or a dim-hom invariant relative to $\X$ is both homological and a dim-hom invariant relative to $\X$.
\end{thm}

\begin{proof}
 Let $M, N \in \mods\Lambda$. It follows from Proposition~\ref{prop:preHom} that if $p(M) = p(N)$ then $q(M) = q(N)$. Thus suppose $q(M) = q(N)$. By Proposition~\ref{prop:homSame} and Proposition~\ref{prop:dimEquiv}, this implies that $[\Hom_\Lambda(T,M)] = [\Hom_\Lambda(T,N)]$ in $K_0(\Gamma)$. Proposition~\ref{cor:grothendieckFree2} then implies that $p(M) = p(N)$ and that $p$ is a homological invariant.
\end{proof}

\begin{rem}\label{rem:preHom}
    This result point towards the possibility that all homological invariants are dim-hom invariants.  Indeed, if we replace the finiteness condition ~\ref{def:homologicalInvariants}(2)(b) in Definition~\ref{def:homologicalInvariants} with the condition that the global dimension of $\Gamma_\X$ is finite, this would follow from Theorem~\ref{thm:hom_dimhom}.  As we suggest in Remark~\ref{rem:exactGroth}(2), we suspect that the two finiteness conditions are equivalent, and thus that all homological invariants are in fact dim-hom invariants.
\end{rem}

%%%%%%%%%%%%%%%%%%%%%%%%%%%

\section{Classes and morphisms of spread modules.}\label{sec:spreads}

We now return to the realm of persistence modules over a finite poset $\poset$, or equivalently to modules over the incidence algebra $\incidence$. We recall the definition and notation of spread modules from Section~\ref{sec:spreadDefinition}. We start by focusing on single-source spread modules and the link between persistence and morphisms from such modules. We then characterize the $\incidence$-linear maps between arbitrary spread modules.

\subsection{Single-source spreads} 

The single-source spread modules defined in Section ~\ref{sec:spreadDefinition} are particularly relevant for persistence theory. To be precise, take an arbitrary module $M \in \mods\incidence$, some vector $v\in M(a)$, and consider the set of vertices where $v$ persists:
$$X(v)=\{x\in \poset | \; M(a,x)(v)\neq 0\}.$$
We observe that $X(v)$ is precisely a single-source spread. The following proposition thus ties the notions of morphism and persistence together.

\begin{prop}\label{prop:spread_of_persistence}
    Let $M \in \mods\incidence$ and let $[a,B] \subseteq \poset$ be a single source spread. Then for $v \in M(a)$, the map $M_{[a,B]}(a)\ni 1_K\mapsto v$ induces a morphism $M_{[a,B]}\to M$ if and only if $[a,B]\supseteq X(v)$. Moreover, the single-source spread module $M_{X(v)}$ is isomorphic to the image of the map $P_a\to M$, induced by $P_a(a)\ni 1_K\mapsto v$.
\end{prop}

\begin{proof}
    If the map $M_{[a,B]}(a)\ni 1_K\mapsto v$ induces a morphism $f:M_{[a,B]}\to M$, then for each vertex $x\in [a,B]$, $f_x\circ M_{[a,B]}(a,x)(1_K)=M(a,x)(v)$. In particular, if $x\in X(v)$, $f_x$ can't be zero, so its domain can't be zero. Then $[a,B]\supseteq X(v)$. Conversely, if this containment is respected, the described map is indeed a morphism.
    
    Now recall that $P_a$, the indecomposable projective module at $a$, is isomorphic to $M_{[a,\infty]}$. Clearly $X(v)\subseteq [a,\infty]$, so the map $P_a(a)\ni 1_K\mapsto v$ induces a morphism $g:P_a\to M$. By definition, $g_x\neq 0$ if and only if $x\in X(v)$. This means $\text{im}(g)$ is a thin module of support $X(v)$. By taking $f_x(1_K)$ as the basis vector for the image at each vertex $x$ we get
    \[ \text{im}(g)(x,y)= \begin{cases}
    1_K & \text{ if } x,y\in X(v)\\
    0 & \text{otherwise}
    \end{cases}
    \]
    So $\text{im}(g)\simeq M_{X(v)}$ as expected.
\end{proof}
  
  Single-source spread modules also admit other characterizations. 
    
\begin{prop}\label{prop:single_source_characterization}
    Let $M\in\mods\incidence$. Then the following are equivalent.
    \begin{enumerate}
        \item $M$ is a single-source spread module.
        \item There exists an indecomposable projective $P \in \mods\incidence$ and a (possibly empty) set $\mathcal{Q}$ of proper projective submodules of $P$ such that
    $M \simeq P/(\sum_{Q \in \mathcal{Q}}Q).$
        \item The projective cover of $M$ is indecomposable.
    \end{enumerate}
\end{prop}

\begin{proof}
    We first note that given $a, b \in \poset$, we have that $P_b \subseteq P_a$ if and only if $a \leq b$.
    
    (1$\implies$2): Write $M = M_{[a,B]}$. Define $$B' = \{x\mid x > a, \ x \not\leq B, \text{ and }\forall y \in \poset: y < x \implies y \leq B\}.$$
    That is, $B'$ consists of those elements of $\poset$ which cover elements of $[a,B]$ but do not themselves lie in $[a,B]$. We conclude that $P_x \subset P_a$ for all $x \in B'$. Moreover, we have $M_{[a,B]} \simeq P_a/(\sum_{x \in B'} P_x)$.
    
    (2$\implies$3): Under the assumption of (2), the projective cover of $M$ is the quotient map $P\rightarrow M$.
    
    (3$\implies$1): Let $p: P \rightarrow M$ be the projective cover of $M$. By assumption, there exists $a \in \poset$ such that $P = P_a$. In particular, we have $\dim P_a(x) \leq 1$ for all $x \in \poset$. Now let $B$ be the set of elements which are maximal in
    $$\{x \in \poset \mid P_a(x) \neq 0 \text{ and } (\ker p)(x) = 0\}.$$
    It is clear that $B$ is a set of incomparable elements of $\poset$ and that $M \simeq M_{[a,B]}$. This completes the proof.
\end{proof}

\begin{rem}
    In the setup of Proposition~\ref{prop:single_source_characterization}(2), those single-source spread modules which satisfy $|\mathcal{Q}| \leq 1$ play a key role in the description of the rank invariant via an exact structure given in \cite{BOO}. We discuss this special case further in Section~\ref{sec:rank}.
\end{rem}

\subsection{Morphisms between spread modules}

We now turn our attention towards computing the Hom-spaces between spread modules. We begin with an example.

\begin{ex}\label{ex:homSpace}
Let $\poset = \{1,2,3,4,5\}\times\{1,2,3\}$ be the $3\times 5$ grid. We will denote elements of $\poset$ using concatenation rather than ordered pairs, e.g. 13 will be used in place of $(1,3)$.
Consider the spreads $[\{13,41\},43]$ (drawn below in dashed red) and $[11,\{23,51\}]$ (drawn below in dotted blue).

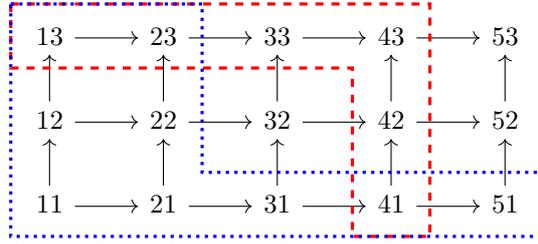
\begin{center}
\begin{tikzcd}[
    every matrix/.style={name=mymatr},
    execute at end picture={
        \draw[very thick,dashed,color=red] ([yshift=5pt,xshift=-5pt]mymatr-1-1.north west)--([yshift=5pt,xshift=5pt]mymatr-1-4.north east)--([yshift=-5pt,xshift=5pt]mymatr-3-4.south east)--([yshift=-5pt,xshift=-5pt]mymatr-3-4.south west)--([yshift=-5pt,xshift=-5pt]mymatr-1-4.south west)--([yshift=-5pt,xshift=-5pt]mymatr-1-1.south west)--cycle;
        \draw[very thick,dotted,color=blue]([yshift=5pt,xshift=-5pt]mymatr-1-1.north west)--([yshift=5pt,xshift=5pt]mymatr-1-2.north east)--([yshift=5pt,xshift=5pt]mymatr-3-2.north east)--([yshift=5pt,xshift=5pt]mymatr-3-5.north east)--([yshift=-5pt,xshift=5pt]mymatr-3-5.south east)--([yshift=-5pt,xshift=-5pt]mymatr-3-1.south west)--cycle;
    }]
    13 \arrow[r] & 23 \arrow[r] & 33 \arrow[r] & 43 \arrow[r] & 53\\
    12 \arrow[r]\arrow[u] & 22 \arrow[r]\arrow[u] & 32 \arrow[r]\arrow[u] & 42 \arrow[r]\arrow[u] & 52\arrow[u]\\
    11\arrow[r]\arrow[u] & 21\arrow[r]\arrow[u] & 31\arrow[r]\arrow[u] & 41\arrow[r]\arrow[u] & 51\arrow[u]
\end{tikzcd}
\end{center}
We want to compute the space of morphisms between the spread representations $M := M_{[\{13,41\},43]}$ and $N := M_{[11,\{23,51\}]}$. Since the supports of $M$ and $N$ only intersect at $13, 23$ and $41$, any $f: M_{[\{13,41\},43]} \rightarrow M_{[11,\{23,51\}]}$ must satisfy $f_v = 0$ for $v \notin \{13,23,41\}$. We then have a pair of commutative diagrams:
\begin{center}
    \begin{tikzcd}
        M &K \arrow[r,"1_K"]\arrow[d,"f_{13}"]& K\arrow[d,"f_{23}"]&K \arrow[r]\arrow[d,"f_{41}"]& 0\arrow[d,"f_{51}"]\\
        N & K \arrow[r,"1_K"]&  K& K \arrow[r,"1_K"]& K.
    \end{tikzcd}
\end{center}
These diagrams imply that $f_{13} = f_{23}$ and that $f_{41} = 0$. Thus choosing a morphism $f:M \rightarrow N$ is equivalent to choosing a scalar $\lambda \in K$ and setting $f_{13}(1) = \lambda$. So, $\text{Hom} \left(M_{[\{13,41\},43]},M_{[11,\{23,51\}]}\right)\simeq K$.
\end{ex}

Generalizing this example, we have the following explicit description of the Hom-space between spread modules. See also \cite[Proposition~3.10]{miller}.

\begin{prop}\label{prop:homSpace}
    Let $[A,B]$ and $[C,D]$ be spreads. Denote by $X_1,\ldots,X_n$ the connected components of $[A,B] \cap [C,D]$ which satisfy
    \begin{equation}\label{eqn:hom_space}\{a\in A \;  | a\leq X_i \} \subseteq X_i \text{ and } \{d \in D \; | \; X_i \leq d\} \subseteq X_i\tag*{$(\star)$}.\end{equation}
    Then there are inverse isomorphisms $\Phi: K^n \rightleftarrows \Hom(M_{[A,B]},M_{[C,D]}): \Psi$
    given as follows.
    \begin{itemize}
        \item Let $(\lambda_1,\ldots,\lambda_n) \in K^n$. We identify each $\lambda_i$ with the linear map $K \rightarrow K$ given as scalar multiplication by $\lambda_i$. Then for each $x \in \poset$, set
        $$\Phi(\lambda_1,\ldots,\lambda_n)_x = \begin{cases} \lambda_i & x \in X_i\\0 & \text{otherwise.}\end{cases}$$
        \item For each $i$, choose some $x_i \in X_i$. Then for $f \in \Hom(M_{[A,B]},M_{[C,D]})$, set
        $$\Psi(f) = (f_{x_1}(1),\ldots,f_{x_n}(1)).$$
    \end{itemize}
\end{prop}

\begin{proof}
    We will first show that $\Psi$ is well defined. Let $f: M_{[A,B]}\rightarrow M_{[C,D]}$, and let
    $x, y \in \poset$ lie in the same connected component $Y$ of $[A,B] \cap [C,D]$. Then $x$ and $y$ are connected by a path in the Hasse quiver of $Y$. Moreover, if $x \leq z \in Y$ then $M_{[A,B]}(x,z) = 1_K = M_{[C,D]}(x,z)$. These two observations together imply that $f_x = f_y$, and therefore $\Psi$ is well defined.
    
    We will next show that $\Psi$ is injective. Let $f: M_{[A,B]}\rightarrow M_{[C,D]}$, and suppose $\Psi(f) = 0$. We wish to show that $f = 0$, or equivalently that $f_x = 0$ for $x \notin \bigcup_{i = 1}^n X_i$. This is clear if $x \notin [A,B]$ or $x \notin [C,D]$, thus assume $x \in [A,B] \cap [C,D]$. By assumption, the connected component $Y$ containing $x$ does not satisfy~\ref{eqn:hom_space}. We then have two cases to consider. Suppose first that there exists $a \in A$ with $a \leq Y$ and $a \notin Y$. Then there exists $y \in Y$ with $a \leq y$. By the previous paragraph, this means $f_a = f_y = f_x$. However, we cannot have $a \in [C,D]$, because otherwise the interval $[a,y]$ would be contained in $Y$. We conclude that $f_x = 0$. The case where there exists $d \in D$ with $Y \leq d$ and $d \notin Y$ is analogous.
    
    We will now show that $\Phi$ is well defined. Let $(\lambda_1,\ldots,\lambda_n) \in K^n$. Obviously, $\Phi(\lambda_1,\ldots,\lambda_n)$ satisfies the morphism condition within each connected component of $[A,B] \cap [C,D]$. The only other morphism conditions to check are at the boundaries of the components $X_1,\ldots,X_n$. Since these components satisfy~\ref{eqn:hom_space}, there is no diagram of either of the following shapes:
    \begin{center}
    \begin{tikzcd}
        K \arrow[r,"1"]\arrow[dd]&K\arrow[dd,"\lambda_i"right] && K\arrow[r]\arrow[dd,"\lambda_i"left]&0\arrow[dd]\\
        \\
        0\arrow[r] & K && K\arrow[r,"1"] & K\\[-1em]
        {\small (\text{vertex outside of $X_i$})} & {\small (\text{vertex inside of $X_i$})} & &{\small (\text{vertex inside of $X_i$})} & {\small (\text{vertex outside of $X_i$})}.
    \end{tikzcd}
\end{center}
Therefore the only patterns that could make the morphism condition fail never happen. We conclude that $\Phi(\lambda_1,\ldots,\lambda_n)$ is well defined.

It is clear that $\Psi \circ \Phi$ is the identity on $K^n$. In particular, this means $\Psi$ is surjective. Since we have already shown that $\Psi$ is injective, we conclude that $\Psi$ is an isomorphism with inverse $\Phi$.
\end{proof}

To conclude this section, we tabulate several consequences of Proposition~\ref{prop:homSpace}.

\begin{coro}\label{cor:bricks}
Let $[A,B]$ be a connected spread. Then $\Hom(M_{[A,B]},M_{[A,B]})= K$. In particular, $M_{[A,B]}$ is indecomposable.
\end{coro}

\begin{proof} The result is a direct consequence of the previous proposition, taking $[A,B] = [C,D]$. Since the intersection of $[A,B]$ with itself is again itself, it is a single connected component. Since $A,B\subseteq [A,B]$, it clearly verifies both conditions of~\ref{eqn:hom_space}.\end{proof}

\begin{coro}\label{coro:homSpace}
    Let $[A,B]$ and $[C,D]$ be spreads and suppose $\Hom(M_{[A,B]},M_{[C,D]}) \neq 0$. Then there exist $a \in A, b \in B, c \in C$, and $d \in D$ such that $c \leq a \leq d \leq b$.
\end{coro}

\begin{proof}
    By Proposition~\ref{prop:homSpace}, there exists a connected component $X$ which satisfies~\ref{eqn:hom_space}. Then there exist $a \in A \cap X$ and $d \in D\cap X$ such that $a \leq d$. Moreover, since $a \in X \subseteq [C,D]$, there exists $c \in C$ with $c \leq a$. Likewise, since $d \in X \subseteq [A,B]$, there exists $b \in B$ with $d \leq b$. This proves the result.
\end{proof}

As a corollary, we also directly get the following characterisation of morphism between intervals, as is well known for type $A$ quivers.

\begin{coro}
There exist a non-zero morphism from $M_{[a,b]}$ to $M_{[c,d]}$ if and only if $c\leq a\leq d\leq b$. In this case, $\emph{Hom}(M_{[a,b]},M_{[c,d]})=K$.
\end{coro}
\begin{proof}
    It follows directly from Corollary ~\ref{coro:homSpace} and Proposition ~\ref{prop:homSpace}.
\end{proof}

Next we show the following necessary conditions to have non-zero morphisms between spreads of the same type, which will be crucial to show Theorem ~\ref{thm:one_source_finite_gldim}. We note that Lemma~\ref{lem:one_source_finite_gldim}(2) below is an immediate consequence of \cite[Proposition~3.10]{miller}, but we provide a proof here for the convenience of the reader.

\begin{lem}\label{lem:one_source_finite_gldim} \
    \begin{enumerate}
        \item Let $[a, B]$ and $[c,D]$ be single-source spreads. Then  $\Hom(M_{[a,B]},M_{[c,D]}) \neq 0$ implies either $c<a$ or $c=a$ and $[c,D]\subseteq [a,B]$.
        \item Let $[A, \infty]$ and $[B,\infty]$ be connected upsets. Then $\Hom (M_{[A,\infty]},M_{[B,\infty]})\neq 0$ implies $[A,\infty]\subseteq [B,\infty]$.    \end{enumerate}
\end{lem}

\begin{proof}
    (1) Suppose $\Hom(M_{[a,B]},M_{[c,D]}) \neq 0$. By Corollary~\ref{coro:homSpace}, we have $c \leq a$. Thus suppose $c = a$. It follows that $X:= [a,B] \cap [c,D]$ is connected and must therefore satisfy~\ref{eqn:hom_space}. It is clear that $c \in X$, so it remains only to show that $D \subseteq X$. Since $c \in X$ and $c \leq d$ for all $d \in D$, we have $\{d\in D| X \leq d\}= D$. As $X$ verifies ~\ref{eqn:hom_space}, we thus have $D\leq X\leq B$. So, $D\leq B$, and $[c,D]\subseteq [a,B]$.
    
    (2) We prove the contrapositive. Suppose $[A,\infty] \not\subseteq [B,\infty]$, and let $f: M_{[A,\infty]} \rightarrow B_{[B,\infty]}$ be a morphism. Because $[A,\infty]$ and $[B,\infty]$ are upsets, the assumption that $[A, \infty] \not \subseteq [B,\infty]$ means there exists $c \in A \setminus [B,\infty]$. In particular, we have $f_c = 0$. Now let $a \in A$ be arbitrary. Since $[A,\infty]$ is connected, it follows that there exists a finite set $\{x_0,\ldots,x_k\} \subseteq [A,\infty]$ such that
    $$c = x_0 \leq x_1 \geq x_2 \leq \cdots \geq x_k = a.$$
    It follows by induction that $f_{x_i} = 0$ for each $x_i$. In particular, we have that $f_a = 0$ for all $a \in A$. By Proposition~\ref{prop:homSpace}, we conclude that $f = 0$.
\end{proof}

Finally, Proposition~\ref{prop:homSpace} can be used to show the following generalisation of the well known fact that irreducible morphisms (in the traditional sense) are either monomorphisms or epimorphisms, see for example \cite[Lemma~5.1]{ARS}. Our proof follows similar arguments to the traditional case, which arise as a special case of the following by setting $\add(X)= \mods\incidence$.

\begin{prop}\label{prop:irreducible}
Suppose that $\add(\X)$ is closed under images in the sense that if $g: R\rightarrow R'$ with $R, R' \in \add(\X)$, then $g(R) \in \add(\X)$.
Let $f: M_{[A,B]} \rightarrow M_{[C,D]}$ be an $\X$-irreducible morphism. Then $f$ is either a monomorphism or an epimorphism. In particular, $\Hom(M_{[A,B]},M_{[C,D]})\simeq K$ and either $[A,B]\subseteq [C,D]$ or $[C,D]\subseteq [A,B]$.
\end{prop}

\begin{proof}
    Consider the factorization $M_{[A,B]} \twoheadrightarrow f(M_{[A,B]}) \hookrightarrow M_{[C,D]}$ of $f$. By assumption, we have $f(M_{[A,B]}) \in \add(\X)$. If the quotient map $M_{[A,B]} \twoheadrightarrow f(M_{[A,B]})$ is a section, then $f$ is a monomorphism. If it isn't, the inclusion map $f(M_{[A,B]}) \hookrightarrow M_{[C,D]}$ must be a retraction. In this case $f$, is an epimorphism.
    
    If $f$ is a monomorphism, then $[A,B] \subseteq [C,D]$ since these are the supports of $M_{[A,B]}$ and $M_{[C,D]}$, respectively. Proposition~\ref{prop:homSpace} then implies that $\Hom(M_{[A,B]},M_{[C,D]}) \simeq K$. The case where $f$ is an epimorphism is completely analogous.
    \end{proof}
    
\begin{rem}
    We note that if we remove the assumption that $\add(\X)$ is closed under images, then there may be $\X$-irreducible morphisms which are neither mono nor epi. For example, in the setting of Example~\ref{ex:homSpace}, take $\X$ to be the set consisting of $M$, $N$, and the indecomposable projectives. Then any non-zero morphism $M \rightarrow N$ is an $\X$-irreducible morphism, but is neither mono nor epi.
\end{rem}

%%%%%%%%%%%%%%%%

%%%%%%%%%%%%%%%%%%%%%%%%%%%%%%%%%%

\section{Homological invariants in persistence theory}\label{sec:newInvariants}
In this section, we use the spread modules discussed in Section~\ref{sec:spreads} to describe examples of homological invariants for persistence modules over some incidence algebra $\incidence$. More precisely, we first use the theory of projectivisation from Section~\ref{sec:projectivisation1} to  provide sufficient conditions for a set $\X$ of spread modules to give rise to a homological invariant. We then use these conditions to prove Theorem~\ref{thm:mainA} as Theorem~\ref{thm:one_source_finite_gldim}  below. In particular, this theorem yields many concrete examples of homological invariants.
Finally, we discuss how our theoretical framework gives an algorithm to compute them. Proposition~\ref{prop:locallyFinite} gives a way to optimize such an algorithm for certain choices of $\X$. It can also be used to extend the reach of this theory to posets that are not finite.

%%%%%%%%%%%%%%%%%%%%%%%%%%%%%%%%%%
\subsection{Homological invariants}\label{sec:newInvariants2}

We now turn our attention to showing that many interesting sets of spread modules give rise to endomorphism algebras of finite global dimension, and thus homological invariants. The basis for these results is the following well-known lemma. Note that this lemma holds more generally, but we have stated it in the generality of the other results of this section. Thus for the duration of this section, $\X$ denotes an arbitrary set of spread modules which contains the indecomposable projectives. 

\begin{lem}\label{lem:acyclic_gldim}
    If the transitive closure of the relation $X \rightarrow Y$ whenever $\Hom_{\incidence}(X,Y)\neq 0$ is a partial order on $\X$, then the endomorphism algebra $\Gamma = \End_{\incidence}(T)^{op}$ has finite global dimension.
\end{lem}

\begin{proof}
    Let $Q$ be the Gabriel quiver of $\Gamma = \End_{\incidence}(T)^{op}$. By the discussion following Definition~\ref{def:irreducible}, there is a directed path from $X$ to $Y$ in $Q$ if and only if there is a sequence of non-zero morphisms
    $$Y = Z_0 \xrightarrow{f_1} Z_1 \xrightarrow{f_2} \cdots \xrightarrow{f_m} Z_m = X$$
    with each $Z_i \in \X$. That is, if $\preceq$ denotes the transitive closure of $\rightarrow$, then there is a directed path from $X$ to $Y$ if and only if $Y \preceq X$. The assumption that $\preceq$ is a partial order then implies that the quiver $Q$ contains no oriented cycles. It is well-known consequence that $\mathsf{gl.dim}\Gamma < \infty$.
\end{proof}

We are now prepared to prove our first main theorem.

\begin{thm}[Theorem~\ref{thm:mainA}]\label{thm:one_source_finite_gldim}
    Let $\poset$ be a finite poset.
    \begin{enumerate}
        \item Let $\X$ be a set of single-source spread modules over $\incidence$ which contains the indecomposable projectives. Then the quotient map $K_0^{split}(\incidence) \rightarrow K_0(\incidence,\X)$ is a homological invariant.
        \item Let $\X$ be a set of connected upset modules over $\incidence$ which contains the indecomposable projectives. Then the quotient map $K_0^{split}(\incidence) \rightarrow K_0(\incidence,\X)$ is a homological invariant.
    \end{enumerate}
\end{thm}

\begin{proof}
       (1) Let $\mathcal{Q}$ be the poset of single-source spreads, ordered by inclusion. By Lemma~\ref{lem:one_source_finite_gldim}(1), the transitive closure of the relation $M_{[a,B]}\rightarrow M_{[c,D]}$ whenever $\Hom(M_{[a,B]},M_{[c,D]})\neq 0$ can be identified with a subposet of the lexicographical product $\mathcal{P}\otimes_{lex}\mathcal{Q}$. The result then follows from Lemma~\ref{lem:acyclic_gldim}  and Proposition~\ref{prop:finiteGlDim}.
       
       (2) Let $\mathcal{U}$ be the set of upsets of $\mathcal{X}$ ordered by inclusion. By Lemma~\ref{lem:one_source_finite_gldim}(2), the transitive closure of the relation $M_{[A,\infty]}\to M_{[B,\infty]}$ whenever $\Hom(M_{[A,\infty]},M_{[B,\infty]})\neq 0$ can be identified with a subposet of $\mathcal{U}$. As with (1), the result then follows from Lemma~\ref{lem:acyclic_gldim}  and Proposition~\ref{prop:finiteGlDim}.
\end{proof}

 \begin{rem}\label{rem:miller}
    The special case of Theorem~\ref{thm:one_source_finite_gldim}(2) where $\X$ is the set of all connected upset modules can also be deduced from \cite[Theorem~6.12]{miller}. Indeed, over a finite poset every (pointwise finite-dimensional) persistence module satisfies \cite[Definition~2.11]{miller}. Thus condition~(6) in \cite[Theorem~6.12]{miller} says precisely that the $\X$-dimension of $\incidence$ is not properly infinite.
\end{rem}

    In case $\X$ is the set of all single-source spread modules, we refer to the quotient map $K_0^{split}(\incidence) \rightarrow K_0(\incidence,\X)$ as the \emph{single-source homological spread invariant}. We defer further discussion of this invariant, and others which fit into the framework of Theorem~\ref{thm:one_source_finite_gldim}(1), to Section~\ref{sec:otherInvariants}. 
    
    It is an interesting problem to characterize precisely which sets of spread modules yield homological invariants. Unfortunately, over arbitrary posets, we cannot hope for an extension of Theorem~\ref{thm:one_source_finite_gldim} to arbitrary choices of $\X$, as the following example shows.

\begin{ex}\label{ex:infinite_gldim}
Let $\mathcal{P}$ be the poset with Hasse diagram
    
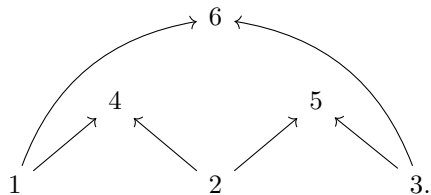
\begin{center}
        \begin{tikzcd}
            &&6\\
            &4&&5\\
            1\arrow[uurr,bend left]\arrow[ur]&&2\arrow[ul]\arrow[ur]&&3.\arrow[uull,bend right]\arrow[ul]
        \end{tikzcd}
    \end{center}
    Note that the incidence algebra $\incidence$ is in this case a path algebra of type $\widetilde{A}_5$. Now set $$\X = \{\text{indecomposable projectives}\} \cup \{M_{[\{1,3\},\{5,6\}]}, M_{[\{1,2\},\{4,6\}]}, M_{[\{2,3\},\{4,5\}]}\}.$$ The Gabriel quiver of $\End_{\incidence}(T)^{op}$ is then
    
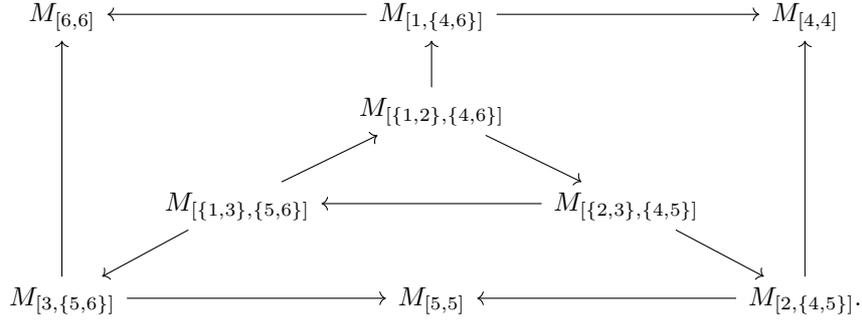
\begin{center}
        \begin{tikzcd}[column sep = 1em]
            M_{[6,6]} \arrow[rr,leftarrow]\arrow[ddd,leftarrow]&&M_{[1,\{4,6\}]}\arrow[d,leftarrow]&&M_{[4,4]}\arrow[ll,leftarrow]\arrow[ddd,leftarrow]\\
            &&M_{[\{1,2\},\{4,6\}]}\arrow[dl,leftarrow]\\
            &M_{[\{1,3\},\{5,6\}]}\arrow[rr,leftarrow]&&M_{[\{2,3\},\{4,5\}]}\arrow[ul,leftarrow]\\
            M_{[3,\{5,6\}]}\arrow[ur,leftarrow]&&M_{[5,5]}\arrow[ll,leftarrow]\arrow[rr,leftarrow]&&M_{[2,\{4,5\}]}.\arrow[ul,leftarrow]
        \end{tikzcd}
    \end{center}
    The oriented 3-cycle at the center of this quiver has the relation that the composition of any two arrows is~0. Moreover, the three oriented pentagons (for example from $M_{\{[1,3\},\{5,6\}]}$ to $M_{[6,6]}$) all have commutativity relations. We claim that the global dimension of $\Gamma$ is infinite. To see this, given a spread module $M_{[A,B]} \in \X$, we denote by $P_{[A,B]}$ the projective $\Gamma$-module at the vertex labeled by $M_{[A,B]}$. Now let $N$ be the indecomposable $\Gamma$-module supported at the vertices labeled by $M_{[\{1,3\},\{5,6\}]}$, $M_{[3,\{5,6\}]}$, and $M_{[5,5]}$. Then $N$ has a minimal projective resolution
        $$\cdots\rightarrow P_{[\{1,2\},\{4,6\}]}\rightarrow P_{[\{1,3\},\{5,6\}]}\rightarrow P_{[\{2,3\},\{4,5\}]}\rightarrow P_{[\{1,2\},\{4,6\}]}\rightarrow N\rightarrow 0.$$
    Indeed, the kernel of the first map is the $\Gamma$-module supported at the vertices labeled by $M_{[\{1,2\},\{4,6\}]}, M_{[1,\{4,6\}]}$, and $M_{[6,6]}$. The symmetry of the quiver then makes it clear that this will be the projective resolution of $N$.
    
    In addition, we can see directly that the $\mathcal{X}$-dimension of $(\mods\incidence,\E_\X)$ is properly infinite. Indeed, the following is readily seen to be a minimal $\X$-resolution of the spread representation $M_{[1,6]}$:
    $$\cdots\rightarrow M_{[\{1,2\},\{4,6\}]}\rightarrow M_{[\{1,3\},\{5,6\}]}\rightarrow M_{[\{2,3\},\{4,5\}]}\rightarrow M_{[\{1,2\},\{4,6\}]}\rightarrow M_{[1,6]}\rightarrow 0.$$
\end{ex}

Note that there are connected spread modules which are not included in set $\X$ in Example~\ref{ex:infinite_gldim}. Thus 
while Example~\ref{ex:infinite_gldim} shows that we cannot fully generalize Theorem~\ref{thm:mainA} to \emph{arbitrary} sets of spread modules, one natural question is what happens when one asks that $\X$ be the set of \emph{all} spread modules.

\begin{openq}\label{conj:gl_dim}\
    Let $\poset$ be an arbitrary finite poset, and let $\mathcal{X}$ be the set of all connected spread modules. Is the quotient map $K_0^{split}(\incidence) \rightarrow K_0(\incidence,\X)$ a homological invariant?
\end{openq}

\begin{rem}\label{rem:repDim}
After a preliminary version of this paper was posted on the arXiv, Asashiba, Escolar, Nakashima, and Yoshiwaki provided a positive answer to Question~\ref{conj:gl_dim} in \cite{AENY2}. Their proof is based upon a result of Iyama \cite{iyama} (see also \cite{ringel}) which shows that for any finite-dimensional algebra $\Lambda$ and any generator-cogenerator $T \in \mods\Lambda$, then there exists $T' \in \mods\Lambda$ such that (i) $\End_\Lambda(T \oplus T')^{op}$ has finite global dimension and (ii) every indecomposable direct summand of $T'$ is a submodule of an indecomposable direct summand of $T$. The positive answer to Question~\ref{conj:gl_dim} then comes from showing that if $\X$ is the set of all connected spreads, then every submodule of an object in $\X$ lies in $\add(\X)$, see \cite[Lemma~4.4]{AENY2}.
\end{rem}
 
%%%%%%%%%%%%%%%%%%%%%%%%%%%%%%
\subsection{Algorithms and extensions to infinite posets}\label{sec:infinite}

Proposition~\ref{prop:finiteGlDim} can be leveraged into an algorithm for computing the value of the invariant $[M]_\X$ when $\Gamma$ has finite global dimension. Indeed, one can first compute a minimal projective resolution of the $\Gamma$-module $\Hom_{\incidence}(T,M)$. By Proposition~\ref{cor:grothendieckFree2}, the class of $M$ in the Grothendieck group can then be computed by alternating sum of the terms of this resolution. The advantage of this approach is that one avoids needing to compute an $\X$-resolution directly in favor of the more straightforward process of computing a projective resolution.

Another way to help compute approximations more effectively is to reduce the size of $\X$ in function of the module we are approximating. In the following proposition, we show how this can be done when all quotients of the modules of $\X$ lie in $\add(\X)$. We emphasize that this is different than saying $\add(\X)$ is closed under quotients. For example, if we take $\X$ to be the set of all single-source spread modules, then all nonzero quotients of the modules in $\X$ will lie in $\X$ (hence also in $\add(\X)$) by Proposition~\ref{prop:single_source_characterization}. On the other hand, the category $\add(\X)$ will not be closed under arbitrary quotients, as this would force $\add(\X) = \mods\incidence$. Indeed, every module is a quotient of its projective cover, which will lie in $\add(\X)$.

\begin{prop}\label{prop:locallyFinite}
 Let $M$ be an arbitrary $\incidence$-module. Let $\X$ a set of spread modules and suppose that every quotient of a module in $\X$ lies in $\add(\X)$. Denote $\Y = \{R \in \X \mid \mathrm{supp}(R) \subseteq \mathrm{supp}(M)\}$. Then every $\X$-approximation
 of $M$ is also a $\Y$-approximation. In particular, the minimal $\X$-approximation of $M$ is exactly the minimal $\Y$-approximation of $M$.
   \end{prop}
 
 \begin{proof}
     It suffices to show that all morphisms from $\add(\X)$ to $M$ factor through morphisms from $\add(\Y)$. Let $R \in \add(\X)$, that is $R=R_1\oplus\cdots\oplus R_n$ for some $R_i\in \X$, and take some morphism $f:R\to M$. Then $f$ admits a factorisation of the form
     $$R\rightarrow \bigoplus_{i = 1}^n f(R_i) \rightarrow \sum_{i = 1}^n f(R_i) = f(R)\hookrightarrow M.$$
     As each $f(R_i)$ is a quotient of a module in $\X$ and satisfies $\mathrm{supp}(f(R_i)) \subseteq \mathrm{supp}(M)$, this shows that $f$ factors through $\add(\Y)$. 
\end{proof}

\begin{rem}
    The previous proposition indicates that much of the theory introduced in this paper may be applicable to finite-dimensional representations of locally finite posets that are not necessarily finite. Note that here we are using finite-dimensional to mean that the total dimension is finite, not to mean pointwise finite-dimensional. In particular, such representations have finite support. Given an infinite set of finite-dimensional spread modules $\X$, we can then use Proposition~\ref{prop:locallyFinite} to form $\X$-approximations of arbitrary finite-dimensional representations. On the other hand, the indecomposable projective representations will in generally not be finite-dimensional, and if $\X$ is infinite then the module $T_\X = \bigoplus_{R \in \X}R$ will not be finite-dimensional either. Thus the extent to which the rest of the paper applies to locally finite posets requires further investigation.
\end{rem}

%%%%%%%%%%%%%%%%%%%%%%%%%%%%%%

\section{Examples and comparison to other invariants}\label{sec:otherInvariants}

In this section, we investigate the relation between our homological invariants and other invariants in persistence theory. It was already shown in Section~\ref{sec:motivation} that the dimension vector is a homological invariant relative to the set of projective modules, see Proposition~\ref{prop:dimEquiv}. We now investigate the barcode, the rank invariant, the generalized rank invariant, the signed barcode, and the generalized persistence diagrams of \cite{KM}.

\subsection{The barcode}\label{sec:barcode}

    In this section, we assume that the Hasse quiver of $\poset = (\poset, \preceq)$ is Dynkin type $A$. More precisely, as a set we identify $\poset$ with $\{1,2,\ldots,|\poset|\}$. We denote by $\preceq$ the partial order defining $\poset$ and by $\leq$ the standard linear order on $\poset = \{1,2,\ldots,|\poset|\}$. The assumption that the Hasse quiver 
    of $\poset$ is Dynkin type $A$ then amounts to the assumption that all cover relations in $\poset$ are either of the form $i \preceq i + 1$ or $i + 1 \preceq i$.
    
    Under this setup, it is well known that the indecomposable $\incidence$-modules are precisely the connected spread modules. Moreover, the connected spreads under the order $\preceq$ correspond directly to the intervals under the order $\leq$. For example, if the quiver of $\poset$ is $1 \rightarrow 2 \leftarrow 3$, it is common to identify the spread module $M_{[\{1,3\},2]}$ with the interval $[1,3]$. More generally, let $F \simeq \mathbb{Z}^{|\poset|\choose 2}$ be the free abelian group generated by the closed intervals of $(\poset, \leq)$. Then the identification between indecomposable $\incidence$-modules and closed intervals of $(\poset, \leq)$ induces an isomorphism $\underline{\mathrm{bar}}:K_0^{split}(\incidence) \rightarrow F$. Given a persistence module $M \in \mods\incidence$, the data $\underline{\mathrm{bar}}[M]$ is often referred to as the \emph{barcode} of $M$.

\begin{prop}\label{prop:barCode}
    Suppose that the Hasse quiver of $\poset$ is Dynkin type $A$. Then the barcode $\displaystyle\underline{\mathrm{bar}}:K_0^{split} \rightarrow \mathbb{Z}^{|\poset|\choose2}$ is a homological invariant.
\end{prop}

\begin{proof}
    Let $\X$ be the set of all connected spread modules. By the above discussion, $\X$ is also the set of all indecomposable $\incidence$-modules. This means the quotient map $K_0^{split}(\incidence)\rightarrow K_0(\incidence,\X)$ is an isomorphism, and is in particular equivalent to $\underline{\mathrm{bar}}$. Thus it remains only to show that the $\X$-dimension of $\incidence$ is finite. To see this, let $T_\X = \bigoplus_{R \in \X}R$. As $\X$ contains all indecomposable modules, $\End_{\incidence}(T_\X)^{op}$ is the so-called \emph{Auslander algebra} of $\incidence$, which is known to have global dimension 2. The result thus follows immediately from Proposition~\ref{cor:grothendieckFree2}.
\end{proof}

%%%%%%%%%%%%%%%%%%%%%%%%%%%%%%%%%%%%%

\subsection{The rank invariant}\label{sec:rank}
We return to $\poset$ denoting an arbitrary finite poset.

Denote by $F$ the free abelian group generated by the (closed) intervals\footnote{Recall that in this paper, intervals are always of the form $[a,b]$ for some $a, b \in \poset$.} of $\poset$. Then the rank invariant (Equation~\ref{eqn:rankInv}) can be considered as a map $p: K_0^{split}(\incidence)\rightarrow F$. It is shown in \cite[Theorem~2.5]{BOO} that $p$ is surjective, and is therefore an invariant in the sense of Definition~\ref{def:invariant}. In order to show that this invariant is homological, we recall a special collection of single-source spreads following \cite{BOO}. We append an extra element $\infty$ to $\poset$, and we set $a < \infty$ for all $a\in \poset$. For any $a < b \in \poset\cup\{\infty\}$, the persistence module $N_{\langle a , b\langle}$, defined by
$$N_{\langle a, b\langle}(c) = \begin{cases} K & a \leq c \text{ and } b \not\leq c\\0 & \text{otherwise}\end{cases}\qquad \qquad N_{\langle a, b\langle}(c,d) = \begin{cases} 1_K & N_{\langle a, b\langle}(c) = K = N_{\langle a, b\langle}(d)\\0 & \text{otherwise,}\end{cases}$$
is called the \textit{hook module} from $a$ to $b$. Note that, if $b \neq \infty$, then $N_{\langle a,b\langle}$ is precisely the quotient of the projective at $a$ by the projective at $b$. Likewise, if $b = \infty$, then $N_{\langle a,b\langle}$ is precisely the projective at $a$. In particular, hook modules are single source spread modules by Proposition~\ref{prop:single_source_characterization}.

One of the main results of \cite{BOO}, specialized to finite posets and restated in the language of this paper, is the following.

\begin{thm}\cite[Theorems~4.4, 4.8 and~4.10]{BOO}\label{thm:boo}
    Let $\poset$ be a finite poset and let $\mathcal{H}$ be the set of hook modules in $\mods\incidence$. Then the rank invariant is  homological with respect to $\mathcal{H}$.
\end{thm}

For the convenience of the reader, we give a proof of this result using the language of our framework.

\begin{proof}
    We first note that the quotient map $K_0^{split}(\incidence) \rightarrow K_0(\incidence,\mathcal{H})$ is homological as a consequence of Theorem~\ref{thm:mainA}(1) and the fact that every hook module is also a single-source spread module. In particular, this means that this quotient map is equivalent to the invariant $p: K_0^{split}(\incidence) \rightarrow \mathbb{Z}^{|\mathcal{H}|}$ given by $p(M) = (\dim\Hom(M,H))_{H \in \mathcal{H}}$. It therefore suffices to show that the rank invariant is equivalent to $p$.
    
    To see this, first let $N_{\langle a,b\langle} \in \mathcal{H}$ with $b \neq \infty$, and let $v$ be a nonzero vector in $N_{\langle a,b\langle}(a)$. Then for any $M \in \mods\incidence$, the morphisms from $N_{\langle a,b\langle}$ to $M$ can be identified with the vectors in $\ker M(a,b)$. Indeed, for $w \in M(a)$ there is a morphism $f: N_{\langle a,b\langle} \rightarrow M$ with $f_a(v) = w$ if and only if $M(a,b)(w) = 0$. As $N_{\langle a,\infty\langle} = P_a$, this means the following all hold:
    \begin{eqnarray*}
        \rk(M(a,b)) &=& \dim\Hom(N_{\langle a,\infty\langle},M) - \dim\Hom(N_{\langle a,b\langle},M)\\
        \dim\Hom(N_{\langle a,b\langle},M) &=& \rk(M(a,a)) - \rk(M(a,b))\\
        \rk(M(a,a)) &=& \dim\Hom(N_{\langle a,\infty\langle},M).
    \end{eqnarray*}
    This shows that $\unrk(M)$ is determined by $p(M)$ and vice versa. Thus $\unrk$ and $p$ are equivalent invariants.
\end{proof}

\begin{rem}
    One consequence of this result, which is made explicit as part of \cite[Theorem~4.4]{BOO}, is that the exact sequences $0 \rightarrow L \rightarrow M \rightarrow N \rightarrow 0$ in $\E_\mathcal{H}$ are characterized by the property that $\unrk(M) = \unrk(L) + \unrk(N)$. Said differently, the set of exact sequences on which $\unrk(-)$ is additive form an exact structure, and this exact structure is precisely $\E_\mathcal{H}$.
\end{rem}

We are now ready to prove our second main theorem, which shows that the single-source spread invariant is finer than the rank invariant.

\begin{thm}[Theorem~\ref{thm:mainC}]\label{thm:finer}\
    Let $\poset$ be a finite poset. Then the single-source homological spread invariant is finer than the rank invariant on $\mods\incidence$. Moreover, these invariants are equivalent if and only if for all $a \in \poset$ the upset $[a,\infty]$ is totally ordered.
\end{thm}

\begin{proof}
    Let $\X$ be the set of single-source spread modules. Since $\X$ contains all of the hook modules in $\mods\incidence$, it follows immediately from Lemma~\ref{lem:finer} and Theorem~\ref{thm:boo} that the single-source spread invariant is finer than the rank invariant.
    
    Now suppose that for all $a \in \poset$, the upset $[a,\infty]$ is totally ordered. We claim that the interval modules, hook modules, and single-source spread modules all coincide. Indeed,
    suppose $[a,B]$ is a single source spread. Then $[a,B] \subseteq [a,\infty]$ is totally ordered, and so $B = \{b\}$ for some $b \in \poset$. In particular, $[a,B] = [a,b]$ is an interval. Now if $[a,b] = [a,\infty]$, then $M_{[a,b]} = P_a$ is a hook module. Otherwise, the set $\mathcal{S} = \{x \in \poset \mid b < x\}$ is nonempty and totally ordered. We can then write $M_{[a,b]} = N_{\langle a,\min(\mathcal{S})\langle}$, which proves the claim. It then follows from Theorem~\ref{thm:boo} that the single-source spread invariant and the rank invariant are equivalent.
    
    Finally, suppose that there exists $a \in \poset$ such that the upset $\mathcal{Q}_a := \{x \in \poset \mid a \leq x\}$ is not totally ordered. Choose two elements $b, c \in \mathcal{Q}_a$ which are incomparable. Without loss of generality, we can assume that $a \leq b$ and $a \leq c$ are both cover relations; that is, we can suppose that $[a,\{b,c\}] = \{a,b,c\}$. Now denote $M := M_{[a,b]} \oplus M_{[a,c]}$ and $M' := M_{[a,a]} \oplus M_{[a,\{b,c\}]}$. As $M$ and $M'$ are not isomorphic and are both direct sums of single-source spread modules, we have $[M]_\X \neq [M']_\X$. On the other hand, we have $\unrk(M) = \unrk(M')$. This concludes the proof.
\end{proof}

\begin{rem}
    The hypotheses of Theorem~\ref{thm:mainC} include the case where $\poset$ is totally ordered, but are slightly more general. Indeed, we only require that for every $a \in \poset$, the upset $[a,\infty]$ is totally ordered. Such subsets are sometimes called ``principle upsets''. For example, consider the poset with Hasse quiver
    $$
    \begin{tikzcd}
        1 \arrow[r] & 2 \arrow[dr]\\
        & 3 \arrow[r] & 4 \arrow[r] & 5.
    \end{tikzcd}
    $$
    The principle upsets of this poset are $\{1,2,4,5\}, \{2,4,5\}, \{3,4,5\}, \{4,5\}$, and $\{5\}$, all of which are totally ordered. Thus the rank invariant and the homological spread invariant coincide for this particular poset.
\end{rem}

In the spirit of the counterexample used to prove Theorem~\ref{thm:finer}, we consider the following.

\begin{ex}\label{ex:recVrk}
    Let $\poset$ be a finite poset which contains a unique maximal element. Let $\mathcal{H}$ be the set of hook modules and let $\mathcal{I}$ be the set of interval modules. Since $\poset$ contains a unique maximal element, we note that $\mathcal{I}$ contains all of the indecomposable projectives. Since every interval is a single-source spread, the quotient map $K_0^{split}(\incidence)\rightarrow K_0(\incidence,\mathcal{I})$ is therefore a homological invariant by Theorem~\ref{thm:one_source_finite_gldim}(1). We also know that the quotient map $K_0^{split}(\incidence)\rightarrow K_0(\incidence,\mathcal{H})$ is equivalent to the rank invariant and is also homological by Theorem~\ref{thm:boo}. 
    
    These two invariants have the same rank (given by the number of pairs $a \leq b \in \poset$), but are not equivalent. For example, let $\poset = \{(0,0), (0,1), (1,0), (1,1)\}$ be the $2\times 2$ grid. We consider the modules
    \begin{center}
        \begin{tikzcd}
            K \arrow[rr] && 0 &&& K \arrow[rr] && 0\\
            K^2 \arrow[u,"\begin{bmatrix} 1 \ 0\end{bmatrix}"]\arrow[rr,"\begin{bmatrix}0 \ 1\end{bmatrix}"]&& K\arrow[u] &&&K^2 \arrow[u,"\begin{bmatrix} 1 \ 0\end{bmatrix}"]\arrow[rr,"\begin{bmatrix}1 \ 0\end{bmatrix}"]&& K.\arrow[u]\\[-1em]
            &M &&&&& M'
        \end{tikzcd}
    \end{center}
    We note that $\undim(M) = \undim(M')$ and that $\unrk(M) = \unrk(M')$. On the other hand, we see that $M$ is a direct sum of interval modules, while $M'$ has an $\mathcal{I}$-resolution
    $$0 \rightarrow M_{[(1,1),(1,1)]} \rightarrow M_{[(0,0),(1,1)]}\oplus M_{[(0,0),(0,0)]} \rightarrow M' \rightarrow 0.$$
    Since this resolution has a non-zero term in degree 1, this means $[M]_\mathcal{I} \neq [M']_\mathcal{I}$.
\end{ex}

We conclude this section with a detailed explanation of why the signed barcode of \cite{BOO} and the homological invariant with respect to $\mathcal{I}$ are not equivalent.

\begin{rem}\label{rem:recVrk}
    Let $\mathcal{H}$, and $\mathcal{I}$ be as in Example~\ref{ex:recVrk}. We denote by $p: K_0^{split}(\incidence) \rightarrow K_0(\incidence,\mathcal{H})$ and $q:K_0^{split}(\incidence) \rightarrow K_0(\incidence,\mathcal{I})$ the two quotient maps. Now recall that the signed barcode of a persistence module $M$ is obtained in \cite{BOO} by applying M\"obius inversion to the classical rank invariant. The result is an invariant which is contained in $K_0(\incidence,\mathcal{I})$ and is equivalent to $p$. More explicitly, there is an isomorphism $\sigma:K_0(\incidence,\mathcal{H}) \rightarrow K_0(\incidence,\mathcal{I})$ such that $\sigma\circ p(M)$ is the signed barcode of $M$. Asking whether the signed barcode and the homological invariant with respect to $\mathcal{I}$ are equivalent is then the same as asking whether $\sigma \circ p = q$. However, if this equality were to hold then we would have that $p$ and $q$ are equivalent, contradicting Example~\ref{ex:recVrk}. In other words, there exists a persistence module $M$ such that $[M]_\mathcal{I}$ and the signed barcode of $M$ are not equivalent.
\end{rem}

%%%%%%%%%%%%%%%%%%%%%%%%%%%%%%%%%%%%%

\subsection{Generalized rank invariants}

We conclude by discussing the generalizations of the rank invariant appearing in \cite{KM,AENY,BOO}. Again $\poset$ denotes an arbitrary finite poset in this section. Given a set $\mathcal{R}$ of connected spread modules, we recall the ``$\mathcal{R}$-rank invariant" $\rk(-,\mathcal{R})$ from Section~\ref{sec:motivation}. We also recall that a ``generalized persistence diagram" $\delta(-,\mathcal{R})$ can be obtained by applying M\"obius inversion to the $\mathcal{R}$-rank invariant. We emphasize that while $\delta(-,\mathcal{R})$ and $\rk(-,\mathcal{R})$ are equivalent as invariants, $\delta(-,\mathcal{R})$ has the advantage of yielding a ``signed approximation'' of a given persistence module by the spreads in $\mathcal{R}$.

To understand the connection between the $\mathcal{R}$-rank invariants, dim-hom invariants, and homological invariants, we first consider the following example.

\begin{prop}\label{prop:gen_rk_grid}
    Let $\poset = \{1,2\} \times \{1,2,3\}$ be the $2\times 3$ grid. Let $\mathcal{R}$ be the set of connected spreads in $\poset$ and let $\X$ be the set of connected spread modules over $\incidence$. Then the generalized rank invariant $\rk(-,\mathcal{R})$ is not a dim-hom invariant relative to $\mathcal{X}$.
\end{prop}

\begin{proof}
    Let $M$ be the following module:

\begin{center}
    \begin{tikzcd}
        K \arrow[r,"1"] & K \arrow[r] & 0\\
        \\
        K \arrow[r,"\text{\scriptsize $\begin{bmatrix}1\\1\end{bmatrix}$}"]\arrow[uu,"1"] & K^2 \arrow[r,"\text{\scriptsize $\begin{bmatrix}0 \ 1\end{bmatrix}$}"]\arrow[uu,"\text{\scriptsize$\begin{bmatrix}1 \ 0\end{bmatrix}$}" right] & K \arrow[uu]
    \end{tikzcd}
\end{center}
Over this poset, indecomposable modules are uniquely determined by their dimension vector. For example, we identify $M$ with the vector $\begin{bmatrix}1 & 1 & 0\\1 & 2 & 1\end{bmatrix}$.

Seeking a contradiction, we compute the ``generalized persistence diagram'' obtained by applying M\"obius inversion to the invariant $\rk(M,\mathcal{R})$. Following \cite[Proposition~3.19]{KM} or \cite[Definition~5.9]{AENY}, this yields
$$\delta(M,\mathcal{R}) = \begin{bmatrix}1&0&0\\1&1&1 \end{bmatrix} + \begin{bmatrix}1&1&0\\1&1&0 \end{bmatrix} + \begin{bmatrix}0&0&0\\0&1&0 \end{bmatrix} - \begin{bmatrix}1&0&0\\1&1&0 \end{bmatrix}.$$
Now denote
$$N := M \oplus \begin{bmatrix}1&0&0\\1&1&0 \end{bmatrix},\qquad\qquad L := \begin{bmatrix}1&0&0\\1&1&1 \end{bmatrix} \oplus \begin{bmatrix}1&1&0\\1&1&0 \end{bmatrix} \oplus \begin{bmatrix}0&0&0\\0&1&0 \end{bmatrix}.$$
Recalling that $\delta(X,\mathcal{R}) = X$ whenever $X$ is a direct sum of spread modules (see \cite[Theorem~
3.14]{KM}), we conclude that $\delta(N,\mathcal{R}) = \delta(L,\mathcal{R})$. On the other hand, denote $X = \begin{bmatrix} 1 & 0 & 0\\1 & 1 & 0\end{bmatrix}$ and observe that
$$
    \dim\Hom_{\incidence}(X,N) = 1, \qquad\qquad
    \dim\Hom_{\incidence}(X,L) = 0.$$
Thus the dim-hom invariant induced by $\X$ can distinguish $N$ and $L$, while the invariant $\delta(-,\mathcal{R})$ cannot. Since the latter is equivalent to the generalized rank invariant $\rk(-,\mathcal{R})$, this proves the result.
\end{proof}

\begin{rem}
    Since this example is small, it is a coincidence that the signed approximation $\delta(M,\mathcal{R})$ preserves the values of $\dim\Hom(-,M_{[A,B]})$ for $[A,B] \in \mathcal{R}$. If one were to replace $M$ with $\begin{bmatrix}1&2&1\\0&1&1\end{bmatrix}$, the result would be an approximation which no longer preserves these values.
\end{rem}

As an immediate consequence, we obtain the following.

\begin{coro}\label{cor:gen_rk_grid}
    Let $\mathcal{P}$ be a poset in which the $2\times 3$ grid is embedded. Let $\mathcal{R}$ be the set of connected spreads in $\poset$ and let $\X$ be the set of connected spread modules over $\incidence$. Then the generalized rank invariant $\rk(-,\mathcal{R})$ is not a dim-hom invariant relative to $\X$.
\end{coro}

\begin{proof}
All of the modules in the proof of Proposition~\ref{prop:gen_rk_grid} embed into persistence modules over the larger poset, and the computation is preserved by this embedding. 
\end{proof}

To summarize, let $\poset, \mathcal{R}$, and $\mathcal{X}$ be as in Corollary~\ref{cor:gen_rk_grid}. We have shown that the generalized rank invariant $\rk(-,\mathcal{R})$, and thus also the signed barcode $\delta(-,\mathcal{R})$, is not a dim-hom invariant \emph{relative to $\mathcal{X}$}. This choice of $\X$ is particularly natural, see Question~\ref{conj:gl_dim}. Indeed, if the quotient map $p: K_0^{split}(\incidence) \rightarrow K_0(\incidence,\X)$ is a homological invariant, then both $p$ and $\rk(-,\mathcal{R})$ assign a persistence module $M$ to a ``signed approximation'' of $M$ by the set of spread modules. Corollary~\ref{cor:gen_rk_grid} would then imply that these approximations do not coincide, even if one allows for a change of basis.

While restricting $\X$ to be composed of only spread modules is a natural choice, it remains an open question whether the generalized rank invariant is either homological or a dim-hom invariant relative to some other set $\mathcal{Y}$. We have shown with a brute-force calculation that this is not the case when $\poset$ is the $2\times 3$ grid, which leads to the following conjecture.

\begin{conj}\label{conj:not_homological}
    Let $\mathcal{P}$ be a finite poset in which the $2\times 3$ grid is embedded, and let $\mathcal{R}$ be the set of connected spreads in $\poset$. Then the generalized rank invariant $\rk(-,\mathcal{R})$ is neither a homological invariant nor a dim-hom invariant (relative to any $\mathcal{Y})$.
\end{conj}

\bibliographystyle{amsalpha}
\bibliography{refs.bib}

\end{document}